\documentclass[final]{siamltex}

\newcommand{\uu}{{\underline u}}
\newcommand{\up}{{\underline p}}
\newcommand{\uq}{{\underline q}}
\newcommand{\uv}{{\underline v}}

\newcommand{\n}{{\bf n}}
\newcommand{\PP}{\mathbf{P}}

\newcommand{\R}{{\mathbb R}}

\usepackage{float}
\usepackage{amsmath}
\usepackage{amssymb}
\usepackage{graphicx} 
\usepackage{psfrag} 
\usepackage{enumerate} 
\usepackage{color}

\def \ds{\displaystyle}
\def \be{\begin{equation}}
\def \ee{\end{equation}}
\def \bea{\begin{eqnarray}}
\def \eea{\end{eqnarray}}
\def \bean{\begin{eqnarray*}}
\def \eean{\end{eqnarray*}}
\newtheorem{theo}{Theorem}
\newtheorem{lem}{Lemma}
\newtheorem{rem}{Remark}

\definecolor{freeblue}{rgb}{0.25,0.41,0.88}

\title{Robin Schwarz algorithm for the NICEM Method: the $\PP_q$ finite element case}

\author{Caroline Japhet\thanks{Universit\'e Paris 13, Sorbonne Paris Cit\'e, LAGA, CNRS UMR 7539,
       99 Avenue J-B Cl\'ement, F-93430 Villetaneuse, France;
  INRIA Paris-Rocquencourt, BP 105, 78153 Le Chesnay, France;
CSCAMM, University of Maryland College Park, MD 20742 USA, E-mail: japhet@math.univ-paris13.fr;
  Partially supported by GNR MoMaS.} \and
  Yvon Maday\thanks{UPMC Univ Paris 06, UMR 7598, Laboratoire Jacques-Louis Lions, F-75005, Paris, France;
  Institut Universitaire de France; and Brown Univ, Division of Applied Maths, Providence, RI, USA,
  E-mail: maday@ann.jussieu.fr. } \and Fr\'ed\'eric
  Nataf\thanks{CNRS, UMR 7598, Laboratoire Jacques-Louis Lions, F-75005, Paris, France;
  UPMC Univ Paris 06, UMR 7598, Laboratoire Jacques-Louis Lions, F-75005, Paris, France,
    E-mail: nataf@ann.jussieu.fr.}}

\begin{document}

\maketitle

\begin{abstract}
In~\cite{GJMN,JMN10} we proposed a new non-conforming domain
decomposition paradigm, the New Interface Cement Equilibrated Mortar
(NICEM) method, based on Schwarz type methods that allows for the use
of Robin interface conditions on non-conforming grids. The error
analysis was done for $\PP_1$ finite elements, in 2D and 3D. In this
  paper, we provide new numerical analysis results that allow to
  extend this error analysis in 2D for piecewise polynomials of higher
  order and also prove the convergence of the iterative algorithm in
  all these cases.
\end{abstract}


\section{Introduction}
The New Interface Cement Equilibrated Mortar (NICEM) method proposed in~\cite{GJMN}
is an equilibrated mortar domain decomposition method
that allows for the use of optimized Schwarz algorithms with Robin
interface conditions on non-conforming grids.
It has been analyzed in~\cite{JMN10} in 2D and 3D for $\PP_1$ elements.

The purpose of this paper is to extend this numerical analysis in 2D
for piecewise polynomials of higher order. We thus establish new numerical analysis results in the frame of finite element approximation and 
also present the iterative algorithm and prove its convergence
in all these cases.

We first consider the problem at the continuous level:\ \ \ Find $u$
such that
\begin{eqnarray}
\label{eq:pbgen}
{\cal L}(u)&=&f \hbox{ in }\Omega\\
\label{eq:pbgen2}
{\cal C}(u)&=&g  \hbox{ on }\partial\Omega
\end{eqnarray}
where ${\cal L}$ and ${\cal C}$ are partial differential equations.
The original Schwarz algorithm is based on a decomposition of the
domain $\Omega$ into overlapping subdomains and the resolution of
Dirichlet boundary value problems in each subdomain. It has been
proposed in~\cite{Lions} to use more general interface/boundary conditions for
the problems on the subdomains in order to use a non-overlapping
decomposition of the domain. The convergence factor is also dramatically
reduced.
More precisely, let $\Omega$ be a ${\cal C}^{1,1}$ (or convex polygon
in 2D or polyhedron in 3D) domain of $I\!\!R^d$, $d=2$ or $3$; we
assume it is decomposed into $K$ non-overlapping subdomains:
$
\overline \Omega = \cup_{k=1}^{K} \overline\Omega^k.
$
We suppose that the subdomains $\Omega^k, \ 1 \le k \le K$ are either
${\cal C}^{1,1}$ or polygons in 2D or polyhedrons in 3D. We
assume also that this decomposition is geometrically conforming in the
sense that the intersection of the closure of two different
subdomains, if not empty, is either a common vertex, a common edge, or
a common face of the subdomains in 3D\footnote[1]{This assumption is not restrictive since in the case of
a partition geometrically non-conforming, the faces can be decomposed in
subfaces to obtain a geometrical conformity }.
Let $\n_k$ be the outward normal from
$\Omega^k$. Let $({\cal B}_{k,\ell})_{1\le k,\ell \le K, k\not= \ell}$
be the chosen transmission conditions on the interface between 
subdomains $\Omega^k$ and $\Omega^\ell$ (e.g. ${\cal B}_{k,\ell}={\partial\ \over \partial
  \n_k}+\alpha_k$).  What we shall call here a Schwarz type method for
the problem (\ref{eq:pbgen})-(\ref{eq:pbgen2}) is its reformulation:
\ \ \ Find $(u_k)_{1\le k\le K}$ such that
\begin{eqnarray}
{\cal L}(u_k)&=&f \hbox{ in }\Omega^k \nonumber\\
{\cal C}(u_k)&=&g  \hbox{ on } \partial\Omega^k \cap \partial\Omega \nonumber\\
{\cal B}_{k,\ell}(u_k)&=&{\cal B}_{k,\ell}(u_\ell)
\hbox{ on }\partial\Omega^k \cap \partial\Omega^\ell, \nonumber
\end{eqnarray}
leading to the iterative procedure
\begin{eqnarray}
{\cal L}(u_k^{n+1})&=&f \hbox{ in }\Omega^k \nonumber\\
{\cal C}(u_k^{n+1})&=&g  \hbox{ on } \partial\Omega^k \cap \partial\Omega \nonumber\\
{\cal B}_{k,\ell}(u_k^{n+1})&=&{\cal B}_{k,\ell}(u_\ell^{n})
\hbox{ on }\partial\Omega^k \cap \partial\Omega^\ell. \nonumber
\end{eqnarray}
The convergence factor of associated Schwarz-type domain decomposition
methods depends largely on the choice of the transmission operators ${\cal B}_{k,\ell}$
(see for instance~\cite{Hagstrom,Nataf.4,GHAD2,GHAD1,Despres.3,Despres.4,Douglas,BenDespres,Widlund,Keyes,Bourdonnaye}
and~\cite{Nataf,Quarteroni}).
More precisely, transmission conditions which
reduce dramatically the convergence factor
of the algorithm have been proposed (see~\cite{Japhet2,Japhet1,JNR}) for a convection-diffusion equation,
where coefficients in second order transmission conditions where optimized. 

On the other hand, the mortar element method,
first introduced in~\cite{BMP}, enables the use of non-conforming
grids, and thus parallel generation of meshes, local adaptive
meshes and fast and independent solvers.
It is also well suited to the use of 
"Dirichlet-Neumann"~(\cite{Quarteroni}), or "Neumann-Neumann"
preconditioned conjugate gradient method applied to the Schur
complement matrix~\cite{Lacour,AMW,TosWid}.  
In~\cite{AJMN}, a new cement to match Robin interface conditions with non-conforming
grids in the case of a finite volume discretization was introduced and analyzed.
Such an approach has been extended
to a finite element discretization in~\cite{GJMN}. A variant has been
independently implemented in~\cite{Vouvakis} for the Maxwell equations,
without numerical analysis. Another approach, in the finite volume case, has been proposed in~\cite{Saas}.

The numerical analysis of the NICEM method proposed in~\cite{GJMN}
is done in~\cite{JMN10} for $\PP_1$ finite elements, in 2D and 3D.
These results are for interface conditions of order 0 (i.e. ${\cal B}_{k,\ell}={\partial\ \over \partial \n_k}+\alpha_k$)
and are the prerequisites for the goal in designing this non-overlapping method for
interface conditions such as Ventcel interface
conditions which greatly enhance the information exchange between subdomains, see~\cite{JMN12} for preliminary results on the extension
of the NICEM method to Ventcel conditions.

The purpose of this paper is first to present a general finite element NICEM method in the case of $\PP_p$ finite elements, with $p\ge 1$ in 2D and $p=1$ in 3D. We also provide a Robin iterative algorithm and prove its convergence. Then,
we present in full details the error analysis in the case of piecewise polynomials of high order in 2D.

In Section~\ref{sec.defmethod}, we describe the NICEM method in 2D and 3D.
Then, in Section~\ref{sec.algo}, we present the iterative algorithm at the continuous and discrete levels,
and we prove, in both cases, the well-posedness and convergence of the iterative method,
for polynomials of low and high order in 2D, and for $\PP_1$ finite elements in 3D.
The convergence is also proven in 3D for $\PP_p$ finite elements, $p\ge 1$,
in a weak sense. In Section~\ref{sec.bestfit2D}  we extend
the error estimates analysis given in~\cite{JMN10} to 2D piecewise polynomials of higher order.
We finally present in Section~\ref{sec:numresults} simulations for two and four
subdomains, that fit the theoretical estimates.

\section{Definition of the method}\label{sec.defmethod}
We consider the following problem : Find $u$ such that
\bea
\label{initial_BVP1}
(Id - \Delta)u &=& f \quad \mbox{in } \Omega \\
\label{initial_BVP2}
              u &=& 0 \quad \mbox{on } \partial{\Omega},
\eea
%
where $f$ is given in $L^2(\Omega)$. \\
The variational statement of the problem
(\ref{initial_BVP1})-(\ref{initial_BVP2})
consists in writing the problem as follows : Find $u \in H^1_0(\Omega)$
such that
\bea\label{initial_VF}
\int_{\Omega} \left(\nabla u \nabla v +uv \right) dx = \int_{\Omega} fvdx,
\quad \forall v \in H^1_0(\Omega).
\eea
%
We introduce the space $H^1_*(\Omega^k)$ defined by
$$H^1_*(\Omega^k) = \{\varphi\in H^1(\Omega^k),\quad
\varphi = 0 \hbox{ over } \partial\Omega\cap\partial
\Omega^k\},$$
and we introduce $\Gamma^{k,\ell}$ the interface of two adjacent
subdomains, $\Gamma^{k,\ell} =
\partial\Omega^k\cap\partial\Omega^\ell.$

It is standard
to note that the space $H^1_0(\Omega)$ can then be identified with the
subspace of the $K$-tuple $\uv=(v_1,...,v_K)$ that are continuous on
the interfaces:
\bea
V = \{\uv=(v_1,...,v_K) \in \prod_{k=1}^K H^1_*(\Omega^k),  \
\forall k,\ell, k \ne \ell, \ 1 \le k,\ell \le K, \ v_k = v_{\ell}
\mbox{ over }
\Gamma^{k,\ell} \}.
\nonumber
\eea
Following~\cite{JMN10}, in order to glue non-conforming grids with Robin transmission
conditions, we impose the constraint $ v_k = v_{\ell}$ over 
$\Gamma^{k,\ell}$
through a Lagrange multiplier in $H^{-1/2}(\partial\Omega^k)$.
The constrained space is then defined as follows
\bea\label{eq:constrainedspace}
    {\cal V} = \ds\lbrace(\uv,\uq)\in
\left(\prod_{k=1}^K H^1_*(\Omega^k)\right)\times \left(\prod_{k=1}^K
H^{-1/2}(\partial\Omega^k)\right), \nonumber\\
\ v_k=v_\ell\hbox{
  and }q_k = - q_\ell \hbox{ over }\Gamma^{k,\ell}, \ \forall k,\ell\rbrace.
\eea
Then, problem (\ref{initial_VF}) is equivalent to the following one (see~\cite{JMN10}):
Find $(\uu,\up) \in {\cal V}$ such that
\begin{eqnarray*}
\label{eq:constraintpb}
\begin{array}{r}
\ds\sum_{k=1}^{K} \int_{\Omega^k} \left( \nabla u_k\nabla v_k +u_kv_k \right) dx
- \sum_{k=1}^{ K}  \
_{H^{-1/2}(\partial\Omega^k)}<p_k,v_k>_{H^{1/2}(\partial\Omega^k)}\hspace{1.2cm}\\
\ds= \sum_{k=1}^{K} \int_{\Omega^k} f_kv_kdx, \quad \forall \uv \in
\prod_{k=1}^KH^1_*(\Omega^k).
\end{array}
\end{eqnarray*}
Being equivalent with the original problem, where $p_k = {\partial
  u\over\partial \n_k}$ over $\partial\Omega^k$, this problem is well
posed. This can also be directly derived from the proof of an inf-sup
condition that follows from the arguments developed hereafter for the
analysis of the iterative procedure.

Note that the Dirichlet-Neumann condition in \eqref{eq:constrainedspace}
is equivalent to the following combined equality
\begin{eqnarray}\label{eq:Robincd}
 p_k + \alpha u_k = - p_{\ell} + \alpha u_{\ell} \quad \mbox{ over }
\Gamma^{k,\ell}, \quad \forall k,\ell.
\end{eqnarray}
As noticed in~\cite{JMN10}, for regular enough function
it is also equivalent to 
\begin{eqnarray}\label{eq:Robincdint}
\hspace{-1cm}
\ \qquad\int_{\Gamma^{k,\ell}}((p_{k}+\alpha u_{k})-(-p_{\ell}+\alpha
u_{\ell})
)\psi_{k,\ell}
= 0,\
\forall \psi_{k,\ell} \in L^2(\Gamma^{k,\ell}), \ \  \forall k,\ell,
\end{eqnarray}
which is the form under which the discrete method is described.

Let us describe the method in the non-conforming discrete case.
\subsection{Discrete case}\label{sec.discretecasepb}
We  introduce now the discrete spaces for piecewise polynomials of higher order
in 2D. Each $\Omega^k$ is provided with
its own mesh ${\cal T}_h^k, \ 1 \le k \le K$, such that
\bea
\overline \Omega^k=\cup_{T \in {\cal T}_h^k} T. \nonumber
\eea
For $T \in {\cal T}_h^k$, let $h_T$ be the diameter of $T$
($h_T=\sup_{x,y \in T} d(x,y)$) and $h$ the discretization parameter
$h=\max_{1 \le k \le K} h_k,$ with $h_k=\max_{T \in {\cal T}_h^k} h_T.$
As noticed in~\cite{JMN10}, for the sake of readability we prefer to use
$h$ instead of $h_k$, but all the analysis could be performed with $h_k$ instead of $h$.
Let $\rho_T$ be the diameter of the  circle (in 2D) or
sphere (in 3D)
inscribed in $T$, then $\sigma_T=\frac{h_T}{\rho_T}$ is a measure of the
non-degeneracy of $T$. We suppose that ${\cal T}_h^k$ is uniformly regular:
there exists $\sigma$ and $\tau$ independent of $h$ such that
$\forall T \in {\cal T}_h^k, \ \sigma_T \le \sigma,$
$\tau h \le h_T .$
We consider that the sets belonging to the meshes are of simplicial type
(triangles), but
the analysis made hereafter can be applied as well for quadrangular meshes.
Let ${\PP}_p(T)$ denote the space of all polynomials defined over $T$
of total degree less than or equal to $p$. 
The finite elements are of Lagrangian type, of class ${\cal C}^0$.
We define over each subdomain two conforming spaces $Y_h^k$ and
$X_h^k$ by:
\bea
Y_h^k&=&\{v_{h,k} \in {\cal C}^0(\overline \Omega^k),
\ \  {v_{h,k}}_{|T} \in {\PP}_p(T), \ \forall T \in {\cal T}_h^k \},
\nonumber\\
X_h^k&=&\{v_{h,k} \in Y_h^k, \  {v_{h,k}}_{|\partial \Omega^k \cap \partial
\Omega}=0\}.\nonumber
\eea
In what follows we assume that the mesh is designed by taking into
account the geometry of the $\Gamma^{k,\ell}$ in the sense that, the
space of traces over each $\Gamma^{k,\ell}$ of elements of $Y_h^k$ is
a finite element space denoted by ${\cal Y}_h^{k,\ell}$.  
Let $k$ be given, the space ${\cal Y}_h^k$ is then the product space of the
${\cal Y}_h^{k,\ell}$ over each $\ell$ such that
$\Gamma^{k,\ell}\not=\emptyset$. With each such interface we associate
a subspace $\tilde W_h^{k,\ell}$ of ${\cal Y}_h^{k,\ell}$ in the same
spirit as in the mortar element method~\cite{BMP} in 2D or~\cite{BBM}
and~\cite{BraessDahmen} in 3D.  To be more
specific, in 2D if the space $X_h^k$ consists of continuous piecewise
polynomials of degree $\le p$, then it is readily noticed that the
restriction of $X_h^k$ to $\Gamma^{k,\ell}$ consists in finite element
functions adapted to the (possibly curved) side $\Gamma^{k,\ell}$ of
piecewise polynomials of degree $\le p$. This side has two end points
that we denote as $x_0^{k,\ell}$ and $x_N^{k,\ell}$ that belong to the
set of vertices of the corresponding triangulation of
$\Gamma^{k,\ell}$ : $x_0^{k,\ell}, x_1^{k,\ell},...,x_{N-1}^{k,\ell},
x_N^{k,\ell}$. The space $\tilde W_h^{k,\ell}$ is then the subspace of
those elements of ${\cal Y}_h^{k,\ell}$ that are polynomials of degree
$\le p-1$ over both $[x_0^{k,\ell}, x_1^{k,\ell}]$ and
$[x_{N-1}^{k,\ell}, x_N^{k,\ell}]$.  As before, the space $\tilde
W_h^{k}$ is the product space of the $\tilde W_h^{k,\ell}$ over each
$\ell$ such that $\Gamma^{k,\ell}\not=\emptyset$.  
Let $\alpha$ be a given positive  real number.
Following~\cite{JMN10}, the discrete constrained space is defined as
\begin{eqnarray}
\label{disc.const}
{\cal V}_h =  \ds\lbrace(\uu_h,\up_h)\in
\left(\prod_{k=1}^K X_h^k\right)\times
\left(\prod_{k=1}^K \tilde W_h^{k}\right), \nonumber\\
\ \qquad\int_{\Gamma^{k,\ell}}((p_{h,k}+\alpha u_{h,k})-(-p_{h,\ell}+\alpha
u_{h,\ell})
)\psi_{h,k,\ell}
= 0,\
\forall \psi_{h,k,\ell} \in \tilde W_h^{k,\ell}
\rbrace,
\end{eqnarray}
and the discrete problem is the following one :
Find $(\uu_h,\up_h) \in {\cal V}_h$ such that\\\\
$\forall \uv_h=(v_{h,1},...v_{h,K}) \in \prod_{k=1}^K X_h^k,$
\bea\label{pbdiscret}
\sum_{k=1}^{K} \int_{\Omega^k} \left( \nabla u_{h,k}\nabla v_{h,k} +u_{h,k}
v_{h,k} \right) dx
- \sum_{k=1}^{K} \int_{\partial\Omega^k} p_{h,k} v_{h,k} ds
= \sum_{k=1}^{K} \int_{\Omega^k} f_k v_{h,k} dx.\hspace{10mm}
\eea
%
The Robin condition \eqref{disc.const} is the discrete counterpart of \eqref{eq:Robincdint}.
\section{Iterative algorithm}\label{sec.algo}
Let us describe the algorithm in the continuous case, and then in the
non conforming discrete case. In both cases, we prove the convergence
of the algorithm towards the solution of the problem.
\subsection{Continuous case}
Let us consider the Robin interface conditions \eqref{eq:Robincd}.
We introduce the following notations: $\ll p,v \gg_{\partial\Omega^k}= _{H^{-1/2}(\partial\Omega^k)}<p,v>_{H^{1/2}(\partial\Omega^k)}$
and $<p,v>_{\Gamma^{k,\ell}}=_{(H_{00}^{1/2}(\Gamma^{k,\ell}))^{\prime}}<p,v>_{H_{00}^{1/2}(\Gamma^{k,\ell})}$.
The algorithm is then defined as follows: let $(u_k^n,p_k^n) \in
H^1_*(\Omega^k) \times H^{-1/2}(\partial\Omega^k)$
be an approximation of $(u,p)$ in $\Omega^k$ at step $n$.
Then, $(u_k^{n+1},p_k^{n+1})$ is the solution in
$H^1_*(\Omega^k) \times H^{-1/2}(\partial\Omega^k)$ of
\bea
\label{algo_continu}
\int_{\Omega^k} \left( \nabla u_k^{n+1}\nabla v_k
+u_k^{n+1}v_k \right) dx
- \ll p_k^{n+1},v_k\gg_{\partial\Omega^k}
= \int_{\Omega^k} f_kv_kdx, \quad \forall v_k \in H^1_*(\Omega^k),  \hspace{8mm}\\
\label{CI_continu}
<p_k^{n+1}+ \alpha u_k^{n+1},v_k>_{\Gamma^{k,\ell}}=
<- p_{\ell}^{n} + \alpha u_{\ell}^{n},v_k>_{\Gamma^{k,\ell}},
\quad \forall v_k \in H_{00}^{1/2}(\Gamma^{k,\ell}). \hspace{8mm}
\eea
%
It is obvious to remark that this series of equations results in uncoupled
problems
set on every $\Omega^k$. Recalling that $f\in L^2(\Omega)$, the strong
formulation is indeed that
\bea
-\Delta u_k^{n+1} + u_k^{n+1} &=& f_k \hspace{1.7cm}\hbox{ over }\Omega^k \nonumber\\
\ds{\partial u_k^{n+1}\over\partial \n_k} + \alpha u_k^{n+1} &=& -p_\ell^n+\alpha
u_\ell^n \quad\hbox{ over } \Gamma^{k,\ell} \nonumber\\
\label{flux_fort}
\ds p_k^{n+1} &=& {\partial u_k^{n+1}\over\partial \n_k} \hspace{1cm}\hbox{ over }
\partial\Omega^k.
\eea
From this strong formulation it is
straightforward to derive by induction that if each
$p^0_k, \ k=1,...,K$, is chosen in $\prod_\ell H^{1/2}(\Gamma^{k,\ell})$,
then,
for each $k$, $1\le k\le K$, and $n\ge 0$ the solution
$u_k^{n+1}$ belongs to
$H^1(\Omega^k)$ and $p_k^{n+1}$
belongs to $\prod_\ell
H^{1/2}(\Gamma^{k,\ell})$ by standard trace results ($p_k^{n+1} = -p_\ell^n+\alpha(u^n_\ell-u_k^{n+1})$). This regularity
assumption on $p^0_k$ will be done hereafter.

  We can prove now that the algorithm (\ref{algo_continu})-(\ref{CI_continu})
converges for all $f \in L^2(\Omega)$:
\begin{theo}
Assume that $f$ is in $L^2(\Omega)$ and $(p^0_k)_{1 \le k \le K}
\in \prod_\ell H^{1/2}(\Gamma^{k,\ell})$. Then, the algorithm
(\ref{algo_continu})-(\ref{CI_continu}) converges
in the sense that
\bea
\lim_{n \longrightarrow \infty} \left( \|u_k^n - u_k\|_{H^1(\Omega^k)}
+  \|p_k^n-p_k\|_{H^{-1/2}(\partial\Omega^k)} \right)=0,
\mbox{  for } 1\le k\le K, \nonumber
\eea
%
where $u_k$ is the restriction to $\Omega^k$ of the solution $u$
to (\ref{initial_BVP1})-(\ref{initial_BVP2}), and $p_k = {\partial
u_k \over\partial \n_k}$ over~$\partial\Omega^k$, $\ 1 \le k \le K$.
\end{theo}
\\\\
{\bf Proof}.
As the equations are linear,
we can take $f=0$. We prove the convergence in the sense that the
associated sequence $(u_k^n,p_k^n)_n$ satisfies
\bea
\lim_{n \longrightarrow \infty} \left( \|u_k^n\|_{H^1(\Omega^k)}
+  \|p_k^n\|_{H^{-1/2}(\partial\Omega^k)} \right)=0,
\mbox{  for } 1\le k\le K. \nonumber
\eea
%
We proceed as in~\cite{Lions,Despres.2} by using an energy estimate that we
derive by taking
$v_k=u_k^{n+1}$ in \eqref{algo_continu} and the use of the regularity
property that
$p_k^{n+1} \in L^2(\partial\Omega^k)$
\bea
\int_{\Omega^k} \left( |\nabla u_k^{n+1}|^2
+|u_k^{n+1}|^2 \right) dx
=  \int_{\partial\Omega^k} p_k^{n+1}u_k^{n+1} ds \nonumber
\eea
%
that can also be written
\bea
\int_{\Omega^k} \left( |\nabla u_k^{n+1}|^2
+|u_k^{n+1}|^2 \right) dx
=\sum_{\ell} \frac{1}{4\alpha} \int_{\Gamma^{k,\ell}}\left(
( p_k^{n+1}+\alpha u_k^{n+1})^2 - ( p_k^{n+1}-\alpha u_k^{n+1})^2\right)ds.
  \nonumber
\eea
%
By using the interface conditions (\ref{CI_continu}) we obtain
\bea\label{estim_en}
\int_{\Omega^k} \left( |\nabla u_k^{n+1}|^2
+|u_k^{n+1}|^2 \right) dx
+\frac{1}{4\alpha}\sum_{\ell}\int_{\Gamma^{k,\ell}}
( p_k^{n+1}-\alpha u_k^{n+1})^2ds
\nonumber\\
= \frac{1}{4\alpha}\sum_{\ell}\int_{\Gamma^{k,\ell}}
( - p_{\ell}^{n}+\alpha u_{\ell}^{n})^2ds.
\eea
%
Let us now introduce two quantities defined at each step $n$ by :
\bea
E^n=\sum_{k=1}^K \int_{\Omega^k} \left( |\nabla u_k^{n}|^2
+|u_k^{n}|^2 \right)
%
\quad
\text{and}
\quad
%
B^n = \frac{1}{4\alpha}\sum_{k=1}^K\sum_{\ell \ne k} \int_{\Gamma^{k,\ell}}
( p_k^{n}-\alpha u_k^{n})^2ds. \nonumber
\eea
%
By summing up the estimates (\ref{estim_en}) over $k=1,...,K$, we have
$E^{n+1} + B^{n+1} \le B^n$,
so~that, by summing up these inequalities, now over $n$, we obtain :
\bea
\sum_{n=1}^{\infty} E^{n} \le B^0. \nonumber
\eea
%
We thus have $\lim_{n \longrightarrow \infty} E^n =0$.
Relation (\ref{flux_fort}) then implies :
\bea
\lim_{n \longrightarrow \infty} \|p_k^n\|_{H^{-1/2}(\partial\Omega^k)}=0,
\mbox{  for } k=1,...,K, \nonumber
\eea
which ends the proof of the convergence of the continuous algorithm.$\qquad \Box$
%
\subsection{Discrete case}\label{sec.discretecase}
%
We first introduce the discrete algorithm defined by: let $(u_{h,k}^n,p_{h,k}^n) \in
X_h^k \times \tilde W_h^{k}$
be a discrete approximation of $(u,p)$ in $\Omega^k$ at step $n$.
Then, $(u_{h,k}^{n+1},p_{h,k}^{n+1})$ is the solution in $X_h^k
\times\tilde W_h^k$ of
\bea
\label{algo_discret}
\int_{\Omega^k} \left( \nabla u_{h,k}^{n+1}\nabla v_{h,k}
+u_{h,k}^{n+1}v_{h,k} \right) dx -
  \int_{\partial\Omega^k}p_{h,k}^{n+1} v_{h,k} ds
= \int_{\Omega^k} f_kv_{h,k}dx  ,\ \forall v_{h,k}\in X_h^k, \hspace{9mm}\\
\label{CI_discret}
\hspace{-1mm}\int_{\Gamma^{k,\ell}} (p_{h,k}^{n+1}+ \alpha u_{h,k}^{n+1})\psi_{h,k,\ell} =
\int_{\Gamma^{k,\ell}} ( -p_{h,\ell}^{n} + \alpha
u_{h,\ell}^{n}) \psi_{h,k,\ell} ,
\quad \forall \psi_{h,k,\ell} \in  \tilde W_h^{k,\ell}.\hspace{8mm}
\eea
In order to analyze the convergence of this iterative scheme, we have to
precise the norms that can be used on the Lagrange multipliers $\up_h$.
For any $\up \in \prod_{k=1}^{K}L^2(\partial \Omega^k)$, in addition to
the natural $L^2$ norm, we can define two better suited norms as follows
\bea
\|\up\|_{-{1\over 2}} = \left(\sum_{k=1}^K
\|p_k\|_{H^{-{1\over 2}}(\partial \Omega^k)}^2 \right)^{1 \over 2} \quad \mbox{and} \quad
\|\up\|_{- {1 \over 2},*} = \left(\sum_{k=1}^K \sum_{\scriptstyle \ell=1
\atop{\atop \scriptstyle \ell \ne k}}^K
\|p_k\|_{H^{-{1\over 2}}_*(\Gamma^{k,\ell})}^2 \right)^{1 \over 2},
\nonumber
\eea
where $\|.\|_{H^{-{1\over 2}}_*(\Gamma^{k,\ell})}$ stands for the dual norm of
${H^{{1\over 2}}_{00}(\Gamma^{k,\ell})}$.
We also need a stability result for the Lagrange multipliers, and refer
to~\cite{BB} in 2D and to~\cite{JMN10} in 3D, in which it is shown that,
\begin{lem}\label{lem.faker}
There exists a constant
$c_*$ such that, for
any $p_{h,k,\ell}$ in $\tilde W_h^{k,\ell}$, there exists an element
$w^{h,k,\ell}$ in
$X_h^k$ that
vanishes  over $\partial\Omega^k\setminus\Gamma^{k,\ell}$ and satisfies
\bea
\label{stab1}
\int_{\Gamma^{k,\ell}} p_{h,k,\ell} w^{h,k,\ell} \ge
\|p_{h,k,\ell}\|^2_{H^{-{1\over
2}}_*(\Gamma^{k,\ell})}
\eea
with a bounded norm
\bea
\label{stab2}
\|w^{h, k,\ell}\|_{H^1(\Omega^k)} \le c_* \|p_{h,k,\ell}\|_{H^{-{1\over
2}}_*(\Gamma^{k,\ell})}.\nonumber
\eea
\end{lem}

Let $\pi_{k,\ell}$ denote the orthogonal projection operator from
$L^2(\Gamma^{k,\ell})$
onto $\tilde W_h^{k,\ell}$. Then, for $v \in L^2(\Gamma^{k,\ell})$,
$\pi_{k,\ell}(v)$ is the unique element of $\tilde W_h^{k,\ell}$
such that
\begin{eqnarray}\label{eq:defpi}
\int_{\Gamma^{k,\ell}} (\pi_{k,\ell}(v)-v)\psi=0, \quad \forall \psi \in
\tilde W_h^{k,\ell}.
\end{eqnarray}

We are now in a position to prove the convergence of the iterative scheme

\begin{theo}\label{theo2}
Let us assume that $\alpha h \le c$, for some small enough constant $c$.
Then, the discrete problem (\ref{pbdiscret}) has a unique solution
$(\uu_h,\up_h) \in {\cal V}_h$.
The algorithm (\ref{algo_discret})-(\ref{CI_discret}) is well posed and
converges
in the sense that
\bea
\lim_{n \longrightarrow \infty} \left( \|u_{h,k}^n - u_{h,k}\|_{H^1(\Omega^k)}
+ \sum_{\ell\neq k} \|p_{h,k,\ell}^n-p_{h,k,\ell}\|_{H^{-{1\over
2}}_*(\Gamma^{k,\ell})} \right)=0,
\mbox{  for } 1\le k\le K. \nonumber
\eea
%
\end{theo}
%
{\bf Proof}. For the sake of convenience, we drop out the index $h$ in what
follows.
We first assume that problems (\ref{pbdiscret}) and
(\ref{algo_discret})-(\ref{CI_discret}) are well posed and
proceed as in the continuous case and assume that $f=0$.
From (\ref{eq:defpi}) we have

$$\forall v_k \in L^2(\Gamma^{k,\ell}), \quad
\int_{\Gamma^{k,\ell}} p_k^{n+1} v_k = \int_{\Gamma^{k,\ell}} p_k^{n+1}
\pi_{k,\ell}(v_k),$$
and (\ref{CI_discret}) also reads
\bea\label{eq:constraintprojn}
p_k^{n+1}+\alpha \pi_{k,\ell} (u_k^{n+1})=
\pi_{k,\ell} (-p_{\ell}^n+\alpha u_{\ell}^n) \quad \mbox{over }
\Gamma^{k,\ell}. 
\eea
%
By taking $v_k=u_k^{n+1}$ in (\ref{algo_discret}), we thus have
\bea
\int_{\Omega^k} \left( |\nabla u_k^{n+1}|^2
+|u_k^{n+1}|^2 \right) dx
\hspace{7cm}\nonumber\\
=  \sum_{\ell} \frac{1}{4\alpha} \int_{\Gamma^{k,\ell}}\left(
( p_k^{n+1}+\alpha \pi_{k,\ell} (u_k^{n+1}))^2 - ( p_k^{n+1}-\alpha
\pi_{k,\ell} (u_k^{n+1}))^2\right)ds.
  \nonumber
\eea
%
Then, by using the interface conditions (\ref{eq:constraintprojn}) we obtain
\bea
\int_{\Omega^k} \left( |\nabla u_k^{n+1}|^2
+|u_k^{n+1}|^2 \right) dx
+\frac{1}{4\alpha}\sum_{\ell}\int_{\Gamma^{k,\ell}}
( p_k^{n+1}-\alpha \pi_{k,\ell}(u_k^{n+1}))^2ds\nonumber\\
= \frac{1}{4\alpha}\sum_{\ell}\int_{\Gamma^{k,\ell}}
  (\pi_{k,\ell}(p_{\ell}^{n}-\alpha u_{\ell}^{n}))^2 ds. \nonumber
\eea
It is straightforward to note that
\bea
\int_{\Gamma^{k,\ell}}
(\pi_{k,\ell}(p_{\ell}^{n}-\alpha u_{\ell}^{n}))^2ds
\le \int_{\Gamma^{k,\ell}}
( p_{\ell}^{n}-\alpha  u_{\ell}^{n})^2 ds\nonumber\\
=  \int_{\Gamma^{k,\ell}}(p_{\ell}^{n} -\alpha \pi_{\ell,k}(u_{\ell}^{n}) +
\alpha \pi_{\ell,k}(u_{\ell}^{n}) -
\alpha  u_{\ell}^{n})^2ds\nonumber\\
= \int_{\Gamma^{k,\ell}}(p_{\ell}^{n} -\alpha \pi_{\ell,k}(u_{\ell}^{n}))^2 +
\alpha ^2(\pi_{\ell,k}(u_{\ell}^{n}) - u_{\ell}^{n})^2ds \nonumber
\eea
since $(Id-\pi_{\ell,k})(u_{\ell}^{n})$ is orthogonal to any element in
$\tilde W_h^{\ell,k}$. For the last term above, we recall that (see~\cite{BMP} in 2D
and~\cite{BBM} or~\cite{BraessDahmen} equation (5.1) in 3D)
\begin{eqnarray*}
\label{eq:propr-pilk}
\int_{\Gamma^{k,\ell}}(\pi_{\ell,k}(u_{\ell}^{n}) - u_{\ell}^{n})^2ds
\le c h \|u_\ell^n\|_{H^{1/2}(\Gamma^{k,\ell})}^2
\le c h \| u_\ell^n\|_{H^1(\Omega^{\ell})}^2.
\end{eqnarray*}
With similar notations as those introduced in the
continuous case, we deduce
\bea
E^{n+1} + B^{n+1} \le c \alpha h E^n + B^n \nonumber
\eea
%
and we conclude as in the continuous case: if $c \alpha h < 1$ then
$\lim_{n\rightarrow\infty}E^n = 0$. The convergence of $u_k^n$ towards
0 in the $H^1$ norm follows. Taking $f=0$ in (\ref{algo_discret}), then using (\ref{stab1})
and the convergence of $u_k^n$ towards 0 in the $H^1$ norm, we derive the convergence of $p_k^n$ in the
$H^{-{1\over 2}}_*(\Gamma^{k,\ell})$ norm.
Note that by having $f=0$ and $(u^n,p^n)=0$ prove that
$(u^{n+1},p^{n+1})=0$ from which we derive that the square problem
(\ref{algo_discret})-(\ref{CI_discret}) is uniquely solvable hence well posed.
Similarly, having $f=0$ and getting rid of the superscripts $n$ and $n+1$
in the previous proof gives (with obvious notations) :
\bea
E + B \le c \alpha h E + B.\nonumber
\eea
%
The existence and uniqueness of a solution of (\ref{pbdiscret}) then results with similar arguments.

In~\cite{JMN10} the well-posedness of (\ref{pbdiscret}) is addressed through a  more direct
proof: let us introduce over
$(\prod_{k=1}^K H^1_*(\Omega^k)\times \prod_{k=1}^K
L^2(\partial\Omega^k))\times \prod_{k=1}^K H^1_*(\Omega^k)$ the
bilinear form
\begin{eqnarray*}
\label{fbs_discret}
\tilde a((\uu,\up), \uv)) = \sum_{k=1}^K
\int_{\Omega^k} \left( \nabla u_k\nabla v_k
+u_k v_k \right) dx -\sum_{k=1}^K
  \int_{\partial\Omega^k}p_k v_k ds.
\end{eqnarray*}
The space $\prod_{k=1}^K H^1_*(\Omega^k)$ is endowed with the norm
\bea
\|\uv\|_* = \left(\sum_{k=1}^K \|v_k\|_{H^1(\Omega^k)}^2 \right)^{1 \over
2}.
\nonumber
\eea
\vspace{1mm}
\begin{lem}\label{lem.infsup}
There exists $c^\prime>0$ and a constant $\beta>0$ such that
\vspace{-3mm}
\begin{eqnarray*}
\label{inf-sup_discret}
\mbox{for } \alpha h\le c^\prime, \quad \forall (\uu_h,\up_h) \in {\cal V}_h ,\ \exists \uv_h\in\prod_{k=1}^K
X_h^k, \hspace{3cm}\nonumber\\
\tilde a((\uu_h,\up_h), \uv_h)) \ge
\beta (\|\uu_h\|_*+ \|\up_h\|_{-{1\over 2},*}) \|\uv_h\|_*.
\end{eqnarray*}
Moreover, we have the continuity argument : there exists a constant
$c>0$ such that
\bea\label{ineq:continuity}
\hspace{-8mm}\forall (\uu_h,\up_h) \in {\cal V}_h ,\ \forall \uv_h\in\prod_{k=1}^K
X_h^k, \quad
\tilde a((\uu_h,\up_h), \uv_h)) \le c (\| \uu_h \|_* +
\|\up_h\|_{-{1 \over 2}}) (\| \uv_h \|_*).
\eea
\end{lem}
This lemma is proven in~\cite{JMN10}, based on Lemma~\ref{lem.faker}.
From Lemma~\ref{lem.infsup}, we have
for any $(\tilde\uu_h,\tilde \up_h)\in{\cal V}_h$,
\bea\label{estimuuh}
\| \uu - \uu_h\|_* +  \|\up - \up_h\|_{-{1 \over 2},*}  \le
c (\| \uu - \tilde\uu_h\|_* +  \|\up - \tilde\up_h\|_{-{1 \over 2}}).
\eea
and we are led to the analysis of the best fit of $(\uu,\up)$
by elements in ${\cal V}_h$.

As noticed in~\cite{JMN10}, it is well known~\cite{BB,BraessDahmen} but unusual that
the inf-sup and continuity conditions involve different norms: the
$\|\cdot\|_{-{1 \over 2}}$ and $\|\cdot\|_{-{1\over 2},*}$ norms. Thus, these two different norms
appear in \eqref{estimuuh} and the best approximation analysis will be done using 
the $\|\cdot\|_{-{1 \over 2}}$ norm, while the error estimates will involve the $\|\cdot\|_{-{1\over 2},*}$ norm.

The analysis of the best fit as been done in~\cite{JMN10} in 2D and 3D for
$\PP_1$ approximations. Let us analyze the best approximation of $(\uu,\up)$ by
elements in ${\cal V}_h$ in the general case of higher order approximations in 2D.

\section{Analysis of the best fit in 2D for higher order approximations}\label{sec.bestfit2D}
In this part we analyze the best approximation of $(\uu,\up)$ by
elements in ${\cal V}_h$.

Following the same lines as in the analysis of the best fit in the $\PP_1$ situation of~\cite{JMN10},
we can prove the following results:
\begin{theo}
\label{best-fit}
Let $u \in H^2(\Omega)\cap H^1_0(\Omega)$, 
be such that $\uu=(u_k)_{1\le k\le K}\in \prod_{k=1}^K H^{2+m}(\Omega^k)$ with  $u_k=u_{|\Omega^k}$,  and
$p-1\ge m \ge 0$. Let us set also
$p_{k,\ell}=\frac{\partial u}{\partial {\bf n}_k}$
over each $\Gamma^{k,\ell}$.
Then there exists $\tilde{\uu}_h$ in
$\prod_{k=1}^K X_h^k$
and $\tilde{\up}_h=(\tilde{p}_{k \ell h})$, with $  \tilde{p}_{k \ell h}
\in \tilde W_h^{k,\ell}$
such that $(\tilde{\uu}_h,\tilde{\up}_h)$ satisfy the coupling condition
(\ref{disc.const}), and
\bea
\| \tilde{\uu}_h -\uu\|_*
&\le& c h^{1+m} \sum_{k=1}^K \| u_k \|_{H^{2+m}(\Omega^k)}
  +{c h^m \over \alpha} \sum_{\ell=1}^K\sum_{k < \ell} \| p_{k,\ell} \|_{H^{{1 \over
2}+m}(\Gamma^{k,\ell})},
\nonumber\\ \nonumber\\
\| \tilde{p}_{k \ell h} - p_{k,\ell} \|_{H^{-{1 \over 2}}(\Gamma^{k,\ell})}
&\le& c\alpha h^{2+m} (\|u_k\|_{H^{2+m}(\Omega^k)}
+\|u_{\ell}\|_{H^{2+m}(\Omega^{\ell})})
\nonumber\\
&& \qquad \qquad  \qquad \qquad\qquad\! + \ c h^{1+m} \| p_{k,\ell} \|_{H^{{1 \over 2}+m}(\Gamma^{k,\ell})}.
\nonumber
\eea
where $c$ is a constant independent of $h$ and $\alpha$.
\end{theo}
%

If we assume more regularity on the normal derivatives on the interfaces,
we have
\begin{theo}
\label{best-fit.2}
Under the assumptions of Theorem \ref{best-fit} and assuming in addition that 
$p_{k,\ell}$ is in $H^{{3 \over 2}+m}(\Gamma_{k,\ell})$.
Then there exists $\tilde{\uu}_h$ in
$\prod_{k=1}^K X_h^k$
and $\tilde{\up}_h=(\tilde{p}_{k \ell h}), \  \tilde{p}_{k \ell h} \in
\tilde W_h^{k,\ell}$
such that $(\tilde{\uu}_h,\tilde{\up}_h)$ satisfy
(\ref{disc.const}), and
%
\begin{multline*}
\| \tilde{\uu}_h -\uu\|_*
\le
c h^{1+m} \sum_{k=1}^K \| u_k \|_{H^{2+m}(\Omega^k)}
+{c h^{m+1} \over \alpha}(\log h)^{\beta(m)} \sum_{\ell=1}^K\sum_{k < \ell} \| p_{k,\ell}
\|_{H^{{3\over 2}+m}(\Gamma^{k,\ell})},
\\
\| \tilde{p}_{k \ell h} - p_{k,\ell} \|_{H^{-{1 \over 2}}(\Gamma^{k,\ell})}
\le
c\alpha h^{2+m} (\|u_k\|_{H^{2+m}(\Omega^k)}
+\|u_{\ell}\|_{H^{2+m}(\Omega^{\ell})}) \\
 + \ c h^{2+m} (\log h)^{\beta(m)}\| p_{k,\ell} \|_{H^{{3 \over2}+m}(\Gamma^{k,\ell})}.
\end{multline*}
where $c$ is a constant independent of $h$ and $\alpha$, and $\beta(m)=0$ if
$m\le p-2$ and $\beta(m)=1$ if $m=p-1$.
\end{theo}
%

The main part of the proof is independent of the degree of the
approximation and is done in~\cite{JMN10}.
Only Lemma~4 in~\cite{JMN10} is dependent of the degree of the approximation and is only proven for
a $\PP_1$ approximation. We prove it for higher order approximations:
\begin{lem}\label{lem_1}
Assume the degree of the finite element approximation $p\le 13$. There
exists two constants $c_1>0$ and $c_2>0$ independent of $h$
such that for all
$\eta_{\ell,k}$ in
${\cal Y}_h^{\ell,k}\cap H_0^1(\Gamma^{k,\ell})$, there exists an element
$\psi_{\ell,k}$ in
$\tilde W_h^{\ell,k}$,   such that
\bea
\label{injectif}
\int_{\Gamma^{k,\ell}}(\eta_{\ell,k}+\pi_{k,\ell}(\eta_{\ell,k}))\psi_{\ell,k}
\ge c_1\|\eta_{\ell,k} \|_{L^2(\Gamma^{k,\ell})}^2,\\
\label{stable}
\| \psi_{\ell,k} \|_{L^2(\Gamma^{k,\ell})} \le
c_2 \| \eta_{\ell,k} \|_{L^2(\Gamma^{k,\ell})}.
\eea
\end{lem}
\begin{rem}
The limit $p\le 13$ is related to the arguments used in the proof we propose for this lemma,
thus, a priori, only technical. We have not found how to alleviate
this limit but actually, for applications, this limit is quite above
what is generally admitted as the optimal range for the degree of the
polynomial in $h-P$ finite element methods. Indeed, as regards the question of
accuracy with respect to run time, the publication \cite{Vos} analyses in
full details\footnote{Of course the answer to that question depends on
the implementation of the discretization method and the exact
properties of the solution to be approximated but this indicates a
tendency that is confirmed by implementation on a large set of other
applications.} and on a variety of problems and regularity of
solutions, the accuracy achieved by low to high order finite element
approximations as a function of the number of degrees of freedom and
of the run time. It appears that the use of degrees between 5 and 8 is
quite competitive which motivates the present analysis.  
\end{rem}

The proof of these results is performed in the following steps.
Note that Lemma~\ref{lem:etapsi_base} below, that generalizes one of the main arguments in the proof of
Lemma~4 in \cite{JMN10} to higher degree in 2d would involve, for a similar generalization in 3d
(see Lemma~7 of \cite{JMN10}), an extension to higher order of the theory developed in \cite{BraessDahmen} that does not
exist yet and goes beyond the scope of the present paper.
\subsection{A first technical result}
\noindent\\
\vspace{-2mm}
\begin{lem}\label{lem:etapsi_base}
Let $1 \le p\le 13$ be an integer. There exists $c$ and $C>0$ such that for all
$\eta\in \PP_p([-1,1])$ s.t.
$\eta(-1)=0$ there exists $\psi\in\PP_{p-1}([-1,1])$ s.t.
$\eta(1)=\psi(1)$, and
\begin{eqnarray*}
J(\psi;\eta):=\int_{-1}^1 (\eta\,\psi -\frac{1}{4}(\eta-\psi)^2)
\ge  c \int_{-1}^1 \eta^2 \quad \text{and} \quad 
\int_{-1}^1 \psi^2 &\le& C \int_{-1}^1 \eta^2.
\end{eqnarray*}
\end{lem}
{\bf Proof}. This lemma has been proven in the case $p=1$ in~\cite{JMN10}.
For $p \ge 2$, we prove it by studying for a given $\eta\in \PP_p([-1,1])$,
$\eta\neq 0$ the maximization problem : \\

Find $\psi\in\PP_{p-1}([-1,1])$ such that
\begin{equation}\label{eq:minconst}
J(\psi;\eta)=\max_{\begin{array}{l}
\varphi \in\PP_{p-1}([-1,1])\\
\varphi(1)=\eta(1)
\end{array}}J(\varphi;\eta).
\end{equation}
The function $J$ is strictly concave in $\varphi$ and there exists a function
satisfying the constraint. This problem admits a solution. The functional
$J(\varphi,\eta)$ being quadratic in $(\varphi,\eta)$ and the constraint being
affine, the optimality condition shows that the problem reduces to a linear
problem the right hand side of which depends linearly of $\eta$. The affine
constraint being of rank one, the problem (\ref{eq:minconst}) admits a
unique solution which depends linearly of $\eta$. Therefore, it makes
sense to introduce the operator:
\[
\begin{array}{rcl}
S: \PP_{p,0}([-1,1]) &\longrightarrow&  \PP_{p-1}([-1,1])\\
\hphantom{S: }\eta &\mapsto& \psi \mathrm{\ solution\ to\
(\ref{eq:minconst})},
\end{array}
\]
where $\PP_{p,0}([-1,1])$ is the set of functions of $\PP_p([-1,1])$ that
vanish at $-1$. In
Lemma~\ref{lem:etapsi_base}, we take $\psi=S(\eta)$.   The operator $S$ is
linear from a finite dimensional space to another so that it is continuous
for any norm on these spaces. Therefore there exists $C>0$ possibly depending on $p$ such that
$\int_{-1}^1
\psi^2 \le C \int_{-1}^1 \eta^2$.
Moreover, the function
\[
\begin{array}{rcl}
H: \PP_{p,0}([-1,1])\backslash \{0\} &\longrightarrow&  \R\\
\hphantom{S: }\eta &\mapsto& \ds\frac{J(S(\eta),\eta)}{\ds\int_{-1}^1
\eta^2}
\end{array}
\]
is continuous and such that $H(\eta)=H(\alpha\eta)$ for any $\alpha\neq
0$. Therefore, it reaches its minimum
which is strictly positive as results from the lemma stated and proven in the next subsection
and the proof of Lemma~\ref{lem:etapsi_base} is complete.$\qquad \Box$
\subsection{Another technical result} 
\label{sec:peq2}
\noindent\\
\vspace{-2mm}
\begin{lem}\label{lem:etapsi}
Let $p\le 13$ and $\eta\in\PP_p([-1,1])$ s.t. $\eta(-1)=0$ and $\eta$ is not the null function.
Then,
$J(S(\eta);\eta) > 0.$ 
\end{lem}

{\bf Proof.} We make use of the Legendre polynomials
\[
L_0(x)=1,\ L_1(x)=x,\ (m+1)L_{m+1}(x)=(2m+1)\,x\,L_m(x) - mL_{m-1}(x),\ m\ge 1.
\]
Let us recall that for any $m\ge 0$,
\[
L_m(1)=1,\ L_m(-1)=(-1)^m, \
\int_{-1}^1 L_m(x)\,L_{m'}(x)\,dx = \delta_{m\,m'}\ds\frac{2}{2m+1}.
\]
The polynomial $\eta$ is decomposed on the Legendre polynomials
\[
\eta=\sum_{m=1}^p \eta_m(L_m+L_{m-1}),
\]
and $\psi=S(\eta)$ is sought in the form
\[
\psi=\sum_{m=0}^{p-1} \psi_m L_m
\]
so that it maximizes the quantity $J(\psi;\eta)$ under the constraint
$\eta(1)=\psi(1)$. This corresponds to the min-max problem
\[
\max_{\psi\in\PP_{p-1}([-1,1])}\min_{\mu\in\R} {\cal L}(\psi,\mu)
\]
where
\[
{\cal L}(\psi,\mu) = J(\psi;\eta)-\mu (\psi(1)-\eta(1)).
\]
We have to prove that the optimal value is positive. The optimality
relations w.r.t
$\psi$ give
\begin{eqnarray*}
\frac{3}{2}(\eta_m+\eta_{m+1})-\frac{1}{2}\psi_m=\mu \frac{2m+1}{2},\
1\le m\le p-1,\quad 
\frac{3}{2}\eta_1-\frac{1}{2}\psi_0=\frac{\mu}{2}.
\end{eqnarray*}
Denoting $\ds R_{p-1}=\sum_{m=0}^{p-1}(2m+1)L_m$, with
$\|R_{p-1}\|_{L^2(]-1,1[)}^2=2p^2$, we get
\begin{equation}\label{eq:psi}
\psi = 3\eta -3\eta_p\,L_p - \mu R_{p-1}.
\end{equation}
Hence, the dual problem writes
\[
\min_{\mu\in\R} G(\mu;\eta),
\]
where
$ G(\mu;\eta) := J(3\eta -3\eta_p\,L_p - \mu
R_{p-1};\eta)-\mu(\psi(1)-\eta(1)) $
and $\psi$ satisfies (\ref{eq:psi}).
After some calculations, $G(\mu;\eta) $ appears a second order polynomial in $\mu$:
\begin{equation}\label{eq:G}
G(\mu;\eta) = \frac{p^2}{2} \mu^2 -\mu (2\eta(1)-3\eta_p)
+(2\|\eta\|_{L^2(]-1,1[)}^2-\frac{9}{2}\frac{\eta_p^2}{2p+1});
\end{equation}
its leading coefficient 
 is positive and its discriminent is proven to be negative in the next lemma, from which we derive that
$\min_\mu G(\mu;\eta)$ is positive and the proof is complete. $\qquad \Box$

\begin{lem}\label{lem:Delta}
For $p\le 13$, the discriminant of (\ref{eq:G}):
\begin{equation*}\label{eq:discr}
\Delta(\eta) := (2\eta(1)-3\eta_p)^2 +
p^2 (-4\|\eta\|_{L^2(]-1,1[)}^2 + 9\frac{\eta_p^2}{2p+1})
\end{equation*}
  is negative if $\eta\in\PP_p([-1,1])$, $\eta(-1)=0$ and $\eta$ is not
the null function.
\end{lem}
\noindent\\[3mm]
{\bf Proof of Lemma~\ref{lem:Delta} in the the case $\pmb{p=2}$}. (the proof  for $3\le p\le 13$,
  is given in Appendix~\ref{appendix:A}).
  
For $p=2$, a direct computation shows that 
\[
\Delta(\eta) =-\frac{80}{3} \eta_1^2-\frac{40}{3} \eta_2\eta_1-\frac{133}{15}\eta_2^2.
\]
The discriminant of the corresponding bilinear form is $-8632/9$. It is negative and the lemma is proven in this case.
$\qquad \Box$

\subsection{Proof of Lemma~\ref{lem_1}}
Using the definition of $\pi_{k,\ell}$, (\ref{eq:defpi}), we derive
\bea
\int_{\Gamma^{k,\ell}}(\eta_{\ell,k}+\pi_{k,\ell}(\eta_{\ell,k}))\psi_{\ell,k}
&=&\int_{\Gamma^{k,\ell}} \eta_{\ell,k}\psi_{\ell,k} + \int_{\Gamma^{k,\ell}}
(\pi_{k,\ell}(\eta_{\ell,k}))^2
+ \int_{\Gamma^{k,\ell}}
\pi_{k,\ell}(\eta_{\ell,k})(\psi_{\ell,k} - \eta_{\ell,k}). \nonumber
\eea
%
Then, using the relation
$\pi_{k,\ell}(\eta_{\ell,k})(\psi_{\ell,k} - \eta_{\ell,k})
\ge - (\pi_{k,\ell}(\eta_{\ell,k}))^2-{1 \over
4}(\psi_{\ell,k} - \eta_{\ell,k})^2$
leads~to
\bea
\int_{\Gamma^{k,\ell}}(\eta_{\ell,k}+\pi_{k,\ell}(\eta_{\ell,k}))\psi_{\ell,k}
\ge \int_{\Gamma^{k,\ell}} \eta_{\ell,k}\psi_{\ell,k}  -{1 \over 4}
\int_{\Gamma^{k,\ell}} (\psi_{\ell,k} - \eta_{\ell,k})^2. \nonumber
\eea
  Remind that we have denoted as
$x_0^{\ell,k},
x_1^{\ell,k},...,x_{N-1}^{\ell,k}, x_N^{\ell,k}$ the vertices of the
triangulation of $\Gamma^{\ell,k}$ that belong to $\Gamma^{\ell,k}$.
By Lemma~\ref{lem:etapsi_base}, and an easy 
scaling argument,
there exists $c, C >0$, $\psi_1 \in \PP_{p-1}([x_0^{\ell,k},x_1^{\ell,k}])$, and
 $\psi_N \in \PP_{p-1}([x_{N-1}^{\ell,k},x_N^{\ell,k}])$, such that
\[
\| \psi_1\|_{L^2(x_0^{\ell,k},x_1^{\ell,k})} + \|
\psi_N\|_{L^2(x_{N-1}^{\ell,k},x_N^{\ell,k})} \le C  (\|
\eta_{\ell,k}\|_{L^2(x_0^{\ell,k},x_1^{\ell,k})} +
\|
\eta_{\ell,k}\|_{L^2(x_{N-1}^{\ell,k},x_N^{\ell,k})}),
\]
$\psi_1(x_1^{\ell,k})=\eta_{\ell,k}(x_1^{\ell,k})$,
$\psi_N(x_{N-1}^{\ell,k})=\eta_{\ell,k}(x_{N-1}^{\ell,k})$ and
\[
\int_{x_0^{\ell,k}}^{x_1^{\ell,k}} (\eta_{\ell,k} \psi_1  -{1 \over 4}
  (\psi_1 - \eta_{\ell,k})^2)
+ \int_{x_{N-1}^{\ell,k}}^{x_N^{\ell,k}} (\eta_{\ell,k} \psi_N  -{1 \over 4}
  (\psi_N - \eta_{\ell,k})^2)
\ge c (\int_{x_0^{\ell,k}}^{x_1^{\ell,k}} \eta_{\ell,k}^2+
\int_{x_{N-1}^{\ell,k}}^{x_N^{\ell,k}}
\eta_{\ell,k}^2 ).
\]
Taking
  $\psi_{\ell,k}$ in
$\tilde W_h^{\ell,k}$ as follows
\bea
\psi_{\ell,k}=\left\{
\begin{array}{l}
  \psi_1
\hbox{ over }
]x_0^{\ell,k}, x_{1}^{\ell,k}[\\
\eta_{\ell,k}\hbox{ over }
]x_1^{\ell,k}, x_{N-1}^{\ell,k}[\\
\psi_N
\hbox{ over }
]x_{N-1}^{\ell,k}, x_{N}^{\ell,k}[\\
\end{array}\right.\nonumber
\eea
proves Lemma~\ref{lem_1} with $c_1=\min(1,c)$ and $c_2=\max(1,C)$.$\qquad \Box$
\subsection{Proof of Theorem~\ref{best-fit}}
We follow the same steps as in the proof of Theorem~2 in~\cite{JMN10}.
Let $u_{kh}^1$ be the unique element of $X_h^k$ defined as follows :
\begin{itemize}
\item
$(u_{kh}^1)_{|\partial \Omega^k}$ is the best fit of $u_k$ over
  $\partial \Omega^k$ in ${\cal Y}_h^{k,\ell}$,
\item
$u_{kh}^1$ at the inner nodes of the triangulation (in $\Omega^k$)
  coincide with the interpolate of $u_k$.
\end{itemize}
Then, it satisfies, for $0\le m\le p-1$,
\bea\label{bestfit_u.2}
\| u_{kh}^1 - u_k\|_{L^2(\partial \Omega^k)}
\le c h^{{3 \over 2}+m} \|u_k\|_{H^{2+m}(\Omega^k)},
\eea
%
from which we deduce that
\bea
\label{bestfit_u}
\| u_{kh}^1 - u_k\|_{L^2(\Omega^k)} + h\| u_{kh}^1 - u_k\|_{H^1(\Omega^k)}
\le c h^{2+m} \|u_k\|_{H^{2+m}(\Omega^k)},
\eea
%
and, from Aubin-Nitsche estimate, 
\bea
\label{bestfit_u.3}
\| u_{kh}^1 - u_k\|_{H^{-{1 \over 2}}(\Gamma^{k,\ell})}
\le c h^{2+m} \|u_k\|_{H^{2+m}(\Omega^k)}.
\eea
%
We introduce  separately the best fit $p_{k \ell h}^1$ of
$p_{k,\ell}=\frac{\partial u}{\partial {\bf n}_k}$
over each
$\Gamma^{k,\ell}$ in $\tilde W_h^{k,\ell}$. Then we have, for $0\le m\le p-1$ 
\bea
\label{bestfit_p}
\| p_{k \ell h}^1 - p_{k,\ell} \|_{L^2(\Gamma^{k,\ell})}
&\le& c h^{{1 \over 2}+m} \| p_{k,\ell} \|_{H^{{1 \over 2}+m}(\Gamma^{k,\ell})},
\\
\label{bestfit_p.2}
\| p_{k \ell h}^1 - p_{k,\ell} \|_{H^{-{1 \over 2}}(\Gamma^{k,\ell})}
&\le& c h^{1+m} \| p_{k,\ell} \|_{H^{{1 \over 2}+m}(\Gamma^{k,\ell})}.
\eea
But there is very few chance that $(\uu_h^1,\up_h^1) \in
\left(\prod_{k=1}^K X_h^k\right)\times \left(\prod_{k=1}^K \tilde
W_h^{k}\right)$ satisfy the
coupling condition (\ref{disc.const}). It misses (\ref{disc.const}) of elements
$\epsilon_{k,\ell}$ and $\eta_{\ell,k}$ such that
\bea
\label{inteps_eta1}
\int_{\Gamma^{k,\ell}}(p_{k\ell h}^1+\epsilon_{k,\ell}+\alpha u_{kh}^1)
\psi_{k,\ell}
&=& \int_{\Gamma^{k,\ell}}(-p_{\ell kh}^1+\alpha \eta_{\ell,k} +\alpha u_{\ell
h}^1) \psi_{k,\ell}
,\ \forall \psi_{k,\ell} \in \tilde W_h^{k,\ell}
\hspace{10mm}\\
\label{inteps_eta2}
\int_{\Gamma^{k,\ell}}(p_{\ell kh}^1+\alpha \eta_{\ell,k} +\alpha u_{\ell
h}^1) \psi_{\ell,k} &=&
\int_{\Gamma^{k,\ell}}(-p_{k\ell h}^1-\epsilon_{k,\ell}+\alpha u_{kh}^1)
\psi_{\ell,k}
,\ \forall \psi_{\ell,k} \in \tilde W_h^{\ell,k}.\hspace{10mm}
\eea
In order to correct that, without polluting
(\ref{bestfit_u.2})-(\ref{bestfit_p.2}), for each couple $(k,\ell)$ we
choose one side, e.g. the smaller indexed one (hereafter we shall 
assume that each couple $(k,\ell)$ is ordered by $k < \ell$).
With this choice, we introduce $\epsilon_{k,\ell} \in
\tilde W_h^{k,\ell}$, $\eta_{\ell,k} \in {\cal Y}_h^{\ell,k} \cap
H_0^1(\Gamma^{k,\ell})$ such that the element $(\tilde{\uu}_h,\tilde{\up}_h)$,
defined by
\bea\label{defptilde}
\tilde{u}_{\ell h}=u^1_{\ell h}+\sum_{k<\ell} {\cal R}_{\ell,k}(\eta_{\ell,k}),
\quad
\tilde{p}_{k \ell h}=p^1_{k \ell h}+ \epsilon_{k,\ell} \quad (\mbox{for } k
< \ell),
\eea
satisfy (\ref{disc.const}). Here ${\cal R}_{\ell,k}$ is a discrete lifting operator
as in~\cite{JMN10} (see also~\cite{Widlund,BG})
that satisfies, with a constant $c$ that is $h$-independent,
that vanishes over
$\partial\Omega^{\ell}\setminus\Gamma^{k,\ell}$ and satisfies
\bea \label{lifting}
\forall w \in {\cal Y}_h^{\ell,k} \cap H_0^1(\Gamma^{k,\ell}),
  \ ({\cal R}_{\ell,k}(w))_{| \Gamma_{k,\ell}}=w, \quad \| {\cal R}_{\ell,k}(w) \|_{H^1(\Omega^{\ell})}
\le c \| w \|_{H^{1 \over 2}_{00}(\Gamma^{k,\ell})}.
\eea
The set of equations
(\ref{inteps_eta1})-(\ref{inteps_eta2}) results in a square system of linear algebraic
equations for $\epsilon_{k,\ell}$ and
$\eta_{\ell,k}$ that can be written as follows
\bea\label{disc.const_2}
\int_{\Gamma^{k,\ell}}(\epsilon_{k,\ell}-\alpha \eta_{\ell,k})\psi_{k,\ell}
&=& \int_{\Gamma^{k,\ell}} e_1 \psi_{k,\ell}
,\ \forall \psi_{k,\ell} \in \tilde W_h^{k,\ell}\\
%
\label{disc.const_3}
\int_{\Gamma^{k,\ell}}(\epsilon_{k,\ell}+\alpha \eta_{\ell,k})\psi_{\ell,k}
&=& \int_{\Gamma^{k,\ell}} e_2 \psi_{\ell,k}
,\ \forall \psi_{\ell,k} \in \tilde W_h^{\ell,k},
\eea
with
\bea
\label{e1-e2}
e_1=-p_{k\ell h}^1-p_{\ell kh}^1+\alpha(u_{\ell h}^1-u_{kh}^1),\quad
e_2=-p_{k\ell h}^1-p_{\ell kh}^1+\alpha(u_{kh}^1-u_{\ell h}^1).\quad
\eea
In~\cite{JMN10}, it is shown that the linear system (\ref{disc.const_2})-(\ref{disc.const_3}) is well posed.

We now estimate $\| \tilde{p}_{k \ell h} - p_{k,\ell} \|_{H^{-{1 \over
2}}(\Gamma^{k,\ell})}$ and $\| \tilde{u}_{\ell h} -
u_{\ell}\|_{H^1(\Omega^{\ell})}$, by first estimating $\|\eta_{\ell,k}
\|_{L^2(\Gamma^{k,\ell})}$:
from (\ref{disc.const_2}) and (\ref{disc.const_3}), we get
\bea\label{epskl-etalk}
\epsilon_{k,\ell}=
\pi_{k,\ell}(\alpha \eta_{\ell,k} +e_1), \quad
\alpha \eta_{\ell,k}=
\pi_{\ell,k}(-\epsilon_{k,\ell} +e_2).
\eea
Injecting the first equation of \eqref{epskl-etalk} in \eqref{disc.const_2}-\eqref{disc.const_3}, we obtain
\bea\label{systeme.eta}
\int_{\Gamma^{k,\ell}}(\eta_{\ell,k}+\pi_{k,\ell}(\eta_{\ell,k}))\psi_{\ell,k}
= {1 \over \alpha} \int_{\Gamma^{k,\ell}} (e_2 -
\pi_{k,\ell}(e_1))\psi_{\ell,k}
,\ \forall \psi_{\ell,k} \in \tilde W_h^{\ell,k}.
\eea
Then, from (\ref{injectif}) and (\ref{systeme.eta}) we get
\bea\label{estim.eta}
c_1\|\eta_{\ell,k} \|_{L^2(\Gamma^{k,\ell})}^2
\le {1 \over \alpha} \|e_2 - \pi_{k,\ell}(e_1)\|_{L^2(\Gamma^{k,\ell})}
\|\psi_{\ell,k}\|_{L^2(\Gamma^{k,\ell})},
\eea
and using (\ref{stable}) in (\ref{estim.eta}) yields
\bea\label{estim.eta2}
\|\eta_{\ell,k} \|_{L^2(\Gamma^{k,\ell})}
\le {c_2 \over \alpha c_1} \|e_2 - \pi_{k,\ell}(e_1)\|_{L^2(\Gamma^{k,\ell})}
\le {c_2 \over \alpha c_1} (\| e_2 \|_{L^2(\Gamma^{k,\ell})}
+\| e_1 \|_{L^2(\Gamma^{k,\ell})}).\qquad \
\eea
Now, from \eqref{e1-e2}, for $i=1,2$
\bea
\|e_i \|_{L^2(\Gamma^{k,\ell})}
\le \|p_{k\ell h}^1+p_{\ell kh}^1\|_{L^2(\Gamma^{k,\ell})}
+\alpha \|u_{\ell h}^1-u_{kh}^1 \|_{L^2(\Gamma^{k,\ell})},
\nonumber
\eea
and recalling that $p_{k,\ell}=\frac{\partial u_k}{\partial {\bf n}_k}=
-\frac{\partial u_\ell}{\partial {\bf n}_{\ell}}=-p_{\ell,k}$ over each
$\Gamma^{k,\ell}$
\bea
\|p_{k\ell h}^1+p_{\ell kh}^1\|_{L^2(\Gamma^{k,\ell})}
&\le& \|p_{k\ell h}^1-p_{k,\ell}\|_{L^2(\Gamma^{k,\ell})}
+\|p_{\ell k h}^1-p_{\ell,k}\|_{L^2(\Gamma^{k,\ell})},
\nonumber\\
\|u_{\ell h}^1-u_{kh}^1 \|_{L^2(\Gamma^{k,\ell})}
&\le& \|u_{k h}^1-u_k \|_{L^2(\Gamma^{k,\ell})}
+\|u_{\ell h}^1-u_{\ell} \|_{L^2(\Gamma^{k,\ell})},
\nonumber
\eea
so that, using (\ref{bestfit_u.2}) and (\ref{bestfit_p}), we derive
for $i=1,2$ and $0\le m\le p-1$
\bea\label{estim.ei}
\| e_i \|_{L^2(\Gamma^{k,\ell})}
\le c \alpha h^{{3 \over 2}+m}
(\|u_k\|_{H^{2+m}(\Omega^k)}+\|u_{\ell}\|_{H^{2+m}(\Omega^{\ell})})
+ch^{{1 \over 2}+m}
\| p_{k,\ell} \|_{H^{{1 \over 2}+m}(\Gamma^{k,\ell})}.\hspace{10mm}
\eea
Thus, (\ref{estim.eta2}) yields, for $0\le m \le p-1$,
\bea\label{estim.eta3}
\|\eta_{\ell,k} \|_{L^2(\Gamma^{k,\ell})}
\le c h^{{3 \over 2}+m}
(\|u_k\|_{H^{2+m}(\Omega^k)}+\|u_{\ell}\|_{H^{2+m}(\Omega^{\ell})})
+{ch^{{1 \over 2}+m} \over \alpha}
\| p_{k,\ell} \|_{H^{{1 \over 2}+m}(\Gamma^{k,\ell})}.\hspace{10mm}
\eea
We can now evaluate $\| \tilde{p}_{k \ell h} - p_{k,\ell} \|_{H^{-{1 \over
2}}(\Gamma^{k,\ell})}$, using the second equation of (\ref{defptilde}) :
\bea\label{eval1}
\| \tilde{p}_{k \ell h} - p_{k,\ell} \|_{H^{-{1 \over 2}}(\Gamma^{k,\ell})}
\le \| \epsilon_{k,\ell} \|_{H^{-{1 \over 2}}(\Gamma^{k,\ell})}
+ \| p_{k \ell h}^1 - p_{k,\ell} \|_{H^{-{1 \over 2}}(\Gamma^{k,\ell})}.
\eea
The term $\| p_{k \ell h}^1 - p_{k,\ell} \|_{H^{-{1 \over
      2}}(\Gamma^{k,\ell})}$ is estimated in (\ref{bestfit_p.2}), so
let us focus on the term $\| \epsilon_{k,\ell} \|_{H^{-{1 \over
      2}}(\Gamma^{k,\ell})}$.  From (\ref{epskl-etalk}) we have,
\bea\label{majoeps}
  \| \epsilon_{k,\ell} \|_{H^{-{1 \over 2}}(\Gamma^{k,\ell})}
\le \alpha\|\eta_{\ell,k} \|_{H^{-{1 \over 2}}(\Gamma^{k,\ell})}
+\|e_1 \|_{H^{-{1 \over 2}}(\Gamma^{k,\ell})}
+ \|(Id-\pi_{k,\ell})(\alpha \eta_{\ell,k} +e_1) \|_{H^{-{1 \over
2}}(\Gamma^{k,\ell})}.\hspace{10mm}
\eea
To evaluate $\|e_1 \|_{H^{-{1 \over 2}}(\Gamma^{k,\ell})}$ we proceed
as for $\|e_1 \|_{L^2(\Gamma^{k,\ell})}$ and from (\ref{bestfit_u.3})
and (\ref{bestfit_p.2}) we have, for $i=1,2$, for $0\le m \le p-1$,
\bea\label{estim.ei.2}
\|e_i \|_{H^{-{1 \over 2}}(\Gamma^{k,\ell})}
\le c \alpha h^{2+m}
(\|u_k\|_{H^{2+m}(\Omega^k)}+\|u_{\ell}\|_{H^{2+m}(\Omega^{\ell})})
+ch^{1+m} \| p_{k,\ell} \|_{H^{{1 \over 2}+m}(\Gamma^{k,\ell})}.\nonumber
\eea
The third term in the right-hand side of (\ref{majoeps}) satisfies
\bea
\|(Id-\pi_{k,\ell})(\alpha \eta_{\ell,k} +e_1) \|_{H^{-{1 \over
2}}(\Gamma^{k,\ell})} \le c \sqrt{h} \|\alpha \eta_{\ell,k} +e_1
\|_{L^2(\Gamma^{k,\ell})}. \nonumber
\eea
Then, using (\ref{estim.eta3}) and (\ref{estim.ei}) yields, for $0\le m \le p-1$,
\bea
\|(Id-\pi_{k,\ell})(\alpha \eta_{\ell,k} +e_1) \|_{H^{-{1 \over
2}}(\Gamma^{k,\ell})}
&\le& c\alpha h^{2+m}
(\|u_k\|_{H^{2+m}(\Omega^k)}+\|u_{\ell}\|_{H^{2+m}(\Omega^{\ell})})\nonumber\\&
&+ \ c h^{1+m} \| p_{k,\ell} \|_{H^{{1 \over 2}+m}(\Gamma^{k,\ell})}.
\nonumber
\eea
In order to estimate the term $\|\eta_{\ell,k} \|_{H^{-{1 \over
2}}(\Gamma^{k,\ell})}$ in (\ref{majoeps}), we
use (\ref{systeme.eta}) and then the symmetry of the operator $\pi_{k,\ell}$:
\bea
2\int_{\Gamma^{k,\ell}} \eta_{\ell,k} \psi_{\ell,k} = \int_{\Gamma^{k,\ell}}
(\psi_{\ell,k}-\pi_{k,\ell}(\psi_{\ell,k}))   \eta_{\ell,k}
+{1 \over \alpha} \int_{\Gamma^{k,\ell}} (e_2 -
\pi_{k,\ell}(e_1)) \psi_{\ell,k}.
\nonumber
\eea
Then, we have
\bea
|\int_{\Gamma^{k,\ell}} \eta_{\ell,k} \psi_{\ell,k}|
\le c\sqrt{h}\|\eta_{\ell,k}\|_{L^2(\Gamma^{k,\ell})}
\|\psi_{\ell,k}\|_{H^{1 \over 2}(\Gamma^{k,\ell})}
+{1\over \alpha}\|e_2 - \pi_{k,\ell}(e_1)\|_{H^{-{1 \over 2}}(\Gamma^{k,\ell})}
\|\psi_{\ell,k}\|_{H^{1 \over 2}(\Gamma^{k,\ell})}
\nonumber
\eea
from which we deduce that
\bea\label{estim.etamoins1demi}
\|\eta_{\ell,k}\|_{H^{-{1 \over 2}}(\Gamma^{k,\ell})}
\le c \sqrt{h}\|\eta_{\ell,k}\|_{L^2(\Gamma^{k,\ell})}+ {c\over \alpha}
\|e_2 - \pi_{k,\ell}(e_1)\|_{H^{-{1
\over 2}}(\Gamma^{k,\ell})}.
\eea
Then, using (\ref{estim.eta3}) and the fact that
\bea\label{estim.e2pie1}
\|e_2 - \pi_{k,\ell}(e_1)\|_{H^{-{1 \over 2}}(\Gamma^{k,\ell})}
\le \|e_2\|_{H^{-{1 \over 2}}(\Gamma^{k,\ell})}
+\|e_1\|_{H^{-{1 \over 2}}(\Gamma^{k,\ell})}
+\|e_1- \pi_{k,\ell}(e_1)\|_{H^{-{1 \over 2}}(\Gamma^{k,\ell})} \nonumber\\
\le \|e_2\|_{H^{-{1 \over 2}}(\Gamma^{k,\ell})}
+\|e_1\|_{H^{-{1 \over 2}}(\Gamma^{k,\ell})}
+c\sqrt{h} \|e_1\|_{L^2(\Gamma^{k,\ell})} \hspace{10mm}
\eea
with (\ref{estim.ei}) and (\ref{estim.ei.2}) yields, for $0\le m \le p-1$,
\bea
\|\eta_{\ell,k}\|_{H^{-{1 \over 2}}(\Gamma^{k,\ell})}
\le c h^{2+m}
(\|u_k\|_{H^{2+m}(\Omega^k)}+\|u_{\ell}\|_{H^{2+m}(\Omega^{\ell})})
+{ch^{1+m} \over \alpha}
\| p_{k,\ell} \|_{H^{{1 \over 2}+m}(\Gamma^{k,\ell})}. \nonumber
\eea
Using the previous inequality in (\ref{majoeps}), (\ref{eval1}) yields, for $0\le m \le p-1$,
\bea\label{estim.p}
\| \tilde{p}_{k \ell h} - p_{k,\ell} \|_{H^{-{1 \over 2}}(\Gamma^{k,\ell})}
\le c\alpha h^{2+m}
(\|u_k\|_{H^{2+m}(\Omega^k)}+\|u_{\ell}\|_{H^{2+m}(\Omega^{\ell})})
+c h^{1+m} \| p_{k,\ell} \|_{H^{{1 \over 2}+m}(\Gamma^{k,\ell})}.\hspace{10mm}
\eea
%
Let us now estimate $\| \tilde{u}_{\ell h} - u_{\ell}\|_{H^1(\Omega^{\ell})}$ :
\bea\label{estim.u2}
\| \tilde{u}_{\ell h} - u_{\ell}\|_{H^1(\Omega^{\ell})}
&\le& \| u_{\ell h}^1 - u_{\ell}\|_{H^1(\Omega^{\ell})}
+ \sum_{k < \ell} \| {\cal R}_{\ell,k}(\eta_{\ell,k}) \|_{H^1(\Omega^{\ell})}
\eea
and from (\ref{lifting}) and an inverse inequality
\bea
\| {\cal R}_{\ell,k}(\eta_{\ell,k}) \|_{H^1(\Omega^{\ell})}
\le c h^{- {1 \over 2}} \|\eta_{\ell,k} \|_{L^2(\Gamma^{k,\ell})}. \nonumber
\eea
Hence, from (\ref{estim.eta3}) we have, for $0\le m \le p-1$,
\bea\label{estim.R}
\| {\cal R}_{\ell,k}(\eta_{\ell,k}) \|_{H^1(\Omega^{\ell})}
\le
ch^{1+m}(\|u_k\|_{H^{2+m}(\Omega^k)}+\|u_{\ell}\|_{H^{2+m}(\Omega^{\ell})})
+{c h^m\over \alpha}  \| p_{k,\ell} \|_{H^{{1 \over 2}+m}(\Gamma^{k,\ell})},
\nonumber
\eea
%
and (\ref{estim.u2}) yields, for $0\le m \le p-1$,
\bea\label{estim.u3}
\| \tilde{u}_{\ell h} - u_{\ell}\|_{H^1(\Omega^{\ell})}
\le ch^{1+m} \|u_{\ell}\|_{H^{2+m}(\Omega^{\ell})}
+ ch^{1+m}\sum_{k < \ell} \|u_k\|_{H^{2+m}(\Omega^k)}
\nonumber\\
+{c h^m \over \alpha} \sum_{k < \ell} \| p_{k,\ell} \|_{H^{{1 \over
2}+m}(\Gamma^{k,\ell})}.
\eea
%
which ends the proof of Theorem~\ref{best-fit}.$\qquad \Box$
\subsection{Proof of Theorem~\ref{best-fit.2}} The proof is the same as
for Theorem~\ref{best-fit}, except that the relation (\ref{bestfit_p})
for $0\le m \le p-1$
is changed in
\bea\label{eq:logh}
\| p_{k \ell h}^1 - p_{k,\ell} \|_{L^2(\Gamma^{k,\ell})}
&\le& c h^{{3 \over 2}+m}\ (\log h)^{\beta(m)} \| p_{k,\ell} \|_{H^{{3 \over
2}+m}(\Gamma^{k,\ell})}.
\eea
%
The proof of \eqref{eq:logh} is given in Appendix~\ref{appendix:B}.
Therefore, (\ref{bestfit_p.2}) is changed in
\bea
\| p_{k \ell h}^1 - p_{k,\ell} \|_{H^{-{1 \over 2}}(\Gamma^{k,\ell})}
&\le& c h^{2+m}\ (\log h)^{\beta(m)}\| p_{k,\ell} \|_{H^{{3 \over
2}+m}(\Gamma^{k,\ell})},
\nonumber
\eea
the inequalities (\ref{estim.p}) and (\ref{estim.u3}) are changed respectively in
\bea
\| \tilde{p}_{k \ell h} - p_{k,\ell} \|_{H^{-{1 \over 2}}(\Gamma^{k,\ell})}
\le c\alpha h^{2+m}
(\|u_k\|_{H^{2+m}(\Omega^k)}+\|u_{\ell}\|_{H^{2+m}(\Omega^{\ell})})\nonumber \\
  +c h^{2+m} \ (\log h)^{\beta(m)} \| p_{k,\ell} \|_{H^{{3\over
2}+m}(\Gamma^{k,\ell})}, \nonumber\\[2mm]
%
\| \tilde{u}_{\ell h} - u_{\ell}\|_{H^1(\Omega^{\ell})}
\le c h^{1+m} \|u_{\ell}\|_{H^{2+m}(\Omega^{\ell})}
+ c h^{1+m} \sum_{k < \ell} \|u_k\|_{H^{2+m}(\Omega^k)}\nonumber \\[-1mm]
+{c h^{1+m} \over \alpha} \ (\log h)^{\beta(m)} \sum_{k < \ell} \| p_{k,\ell}
\|_{H^{{3 \over 2}+m}(\Gamma^{k,\ell})}. \qquad \Box
\nonumber
\eea
%
%
\subsection{Error Estimates}
Thanks to (\ref{estimuuh}), we have the following error estimates:
\begin{theo}
\label{error-estimate}
Assume that the solution $u$ of
(\ref{initial_BVP1})-(\ref{initial_BVP2}) is in $H^2(\Omega)\cap H^1_0(\Omega)$, and $u_k=u_{|\Omega^k}\in
H^{2+m}(\Omega^k)$, with
$p-1\ge m \ge 0$,
and let
$p_{k,\ell}=\frac{\partial u}{\partial {\bf n}_k}$
over each $\Gamma^{k,\ell}$.
Then, there exists a constant $c$ independent of $h$ and $\alpha$
such that
\bea
\| \uu_h -\uu\|_* + \| \up_h - \up \|_{-{1 \over 2},*}
\le c(\alpha h^{2+m}+h^{1+m} )
\sum_{k=1}^K \| \uu\|_{H^{2+m}(\Omega^k)} \nonumber\\[-1mm]
+ \ c ({h^m\over \alpha} + h^{1+m})\sum_{k=1}^K \sum_{\ell}\| 
p_{k,\ell} \|_{H^{{1
\over 2}+m}(\Gamma^{k,\ell})}.
\nonumber
\eea
%
\end{theo}
%
\begin{theo}
\label{error-estimate2}
Assume that the solution $u$ of
(\ref{initial_BVP1})-(\ref{initial_BVP2}) is in $H^2(\Omega)\cap H^1_0(\Omega)$, $u_k=u_{|\Omega^k}\in
H^{2+m}(\Omega^k)$,
and
$p_{k,\ell}=\frac{\partial u}{\partial {\bf n}_k}$
is in $H^{{3 \over 2}+m}(\Gamma_{k,\ell})$ with $p-1\ge m \ge 0$.
Then there exists a constant $c$ independent of $h$ and $\alpha$
such that
\bea
\| \uu_h -\uu\|_* + \| \up_h - \up\|_{-{1 \over 2},*}
\le c(\alpha h^{2+m}+h^{1+m} )
\sum_{k=1}^K \| \uu\|_{H^{2+m}(\Omega^k)} \nonumber\\
+ \ c ({h^{1+m}\over \alpha} + h^{2+m}) (\log h)^{\beta(m)} 
\sum_{k=1}^K \sum_{\ell}\|
p_{k,\ell}
\|_{H^{{3
\over 2}+m}(\Gamma^{k,\ell})}
\nonumber
\eea
with $\beta(m)=0$ if
$m\le p-2$ and $\beta(m)=1$ if $m=p-1$.
\end{theo}
%
%
\section{Numerical results}\label{sec:numresults} 
%
We consider the initial problem, with exact solution $u(x,y)=x^4y^4+xy\cos(10xy)$.
The domain is the unit square $\Omega=(0,1) \times (0,1)$,
decomposed into non-overlapping subdomains with meshes generated independently.
In Sections~\ref{subsec:errP2}, \ref{subsec:errP123}, to observe the numerical error estimates for the 
discrete problem (\ref{pbdiscret}), one need to compute the converged solution of the discrete algorithm
(\ref{algo_discret})-(\ref{CI_discret}) regardless of the algorithm used to compute it. Thus it is
the solution at convergence of the algorithm (\ref{algo_discret})-(\ref{CI_discret}) with a stopping criterion on
the residual (i.e. the jumps of interface conditions) that must be extremely small, e.g. smaller than $10^{-14}$.
For all the other simulations where we are interested in $u_h^n$ and not $u_h$, a residual of $10^{-2}$ times the
target $H^1$ error is considered, for stopping the iterations.
%
\subsection{Choice of the Robin parameter}\label{Sec:alpha}
In our simulations the Robin parameter $\alpha$ is either an arbitrary constant
or is obtained by minimizing the convergence factor (depending on the
mesh size in that case, see~\cite{JMN10}).  In the conforming two subdomains case, with
constant mesh size $h$ and an interface of length $L$, the optimal theoretical value of $\alpha$
which minimizes the convergence factor at the continuous level is (see~\cite{Gander06}):
\vspace{-5mm}
\begin{equation}\label{eq.alphaopt}
\alpha_{\text{opt}}(L,h)=[((\frac{\pi}{L})^2+1)((\frac{\pi}{h})^2+1)]^{\frac{1}{4}}.
\end{equation}
Note that this optimal choice for $\alpha$ does not seem to provide an
optimal error estimate from Theorem \ref{error-estimate}.
Nevertheless, as was illustrated in \cite{JMN10}
the regularity of the normal derivative of $u$ along the interfaces
enters most of the times in the frame of Theorem~\ref{error-estimate2}
that allows a larger range of choice for $\alpha$, compatible with the
above mentioned optimal choice (as regards the algorithm).

In the non-conforming case, we consider the following values : 
$\alpha_{\text{min}}=\alpha_{\text{opt}}(L,{h_{\text{min}} \over p})$,
$\alpha_{\text{mean}}=\alpha_{\text{opt}}(L,{h_{\text{mean}}\over p})$,
$\alpha_{\text{max}}=\alpha_{\text{opt}}(L,{h_{\text{max}}\over p})$,
where $h_{\text{min}}$, $h_{\text{mean}}$ and $h_{\text{max}}$ stands respectively for the
smallest, meanest or highest step size on the interface and $p$ is the degree of the approximation.

\subsection{$H^1$ error between the continuous and discrete solutions for $\PP_2$ finite elements}
\label{subsec:errP2}
In this part, we compare the relative $H^1$ error in the
non-conforming case to the error obtained on a uniform conforming
grid.

We define the relative $H^1$ error as follows:
Let $u_k=u_{|\Omega^k}, \ 1 \le k \le K$
(where u is the continuous solution), and let
$(\uu_h)_k=(\uu_h)_{|\Omega^k}$ where $\uu_h$ is the solution of the
discrete problem (\ref{pbdiscret}).  Now, let $N_x=\|u\|_*$ and let
$E_k=\|(\uu_h)_k-u_k\|_{H^1(\Omega^k)}$, $1 \le k \le K$. Let
$E=(\sum_{i=1}^K E_i^2)^{1/2}.$ The relative $H^1$ error is then
$E/N_x$.

\begin{figure}[H]
  \centering 
  \includegraphics[height=3.85cm]{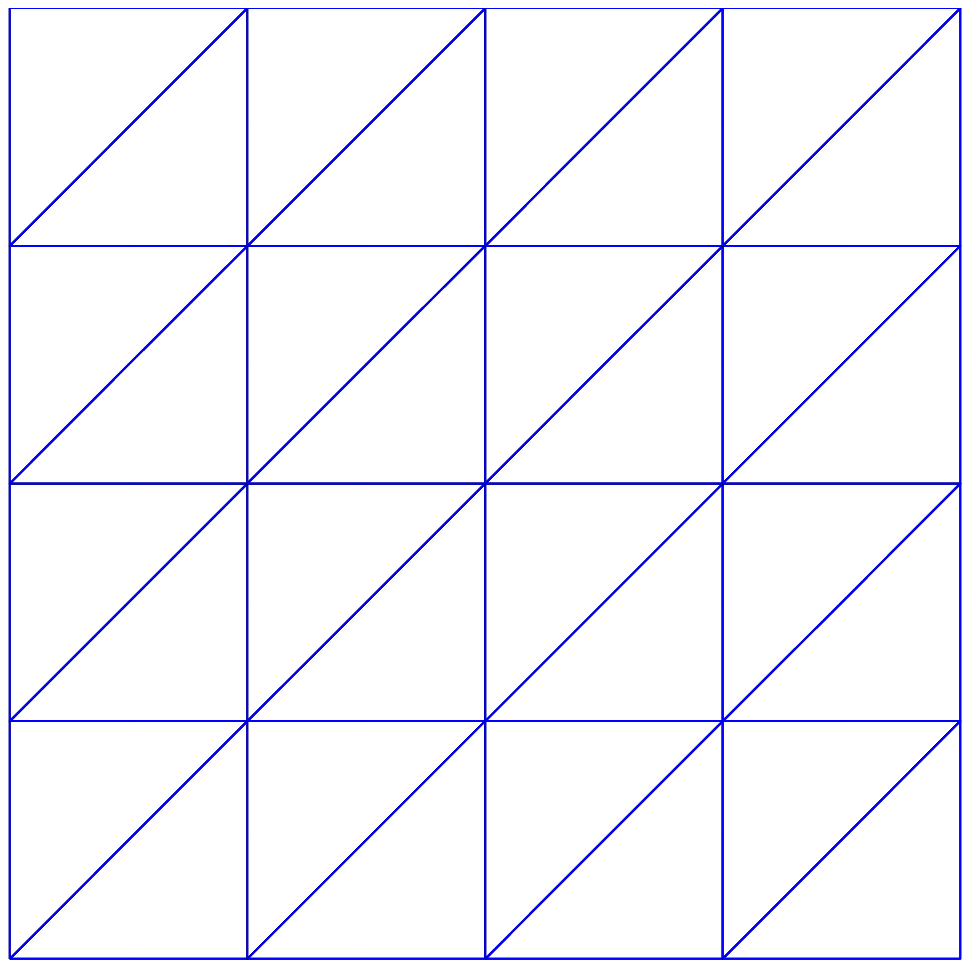}
  \includegraphics[height=3.85cm]{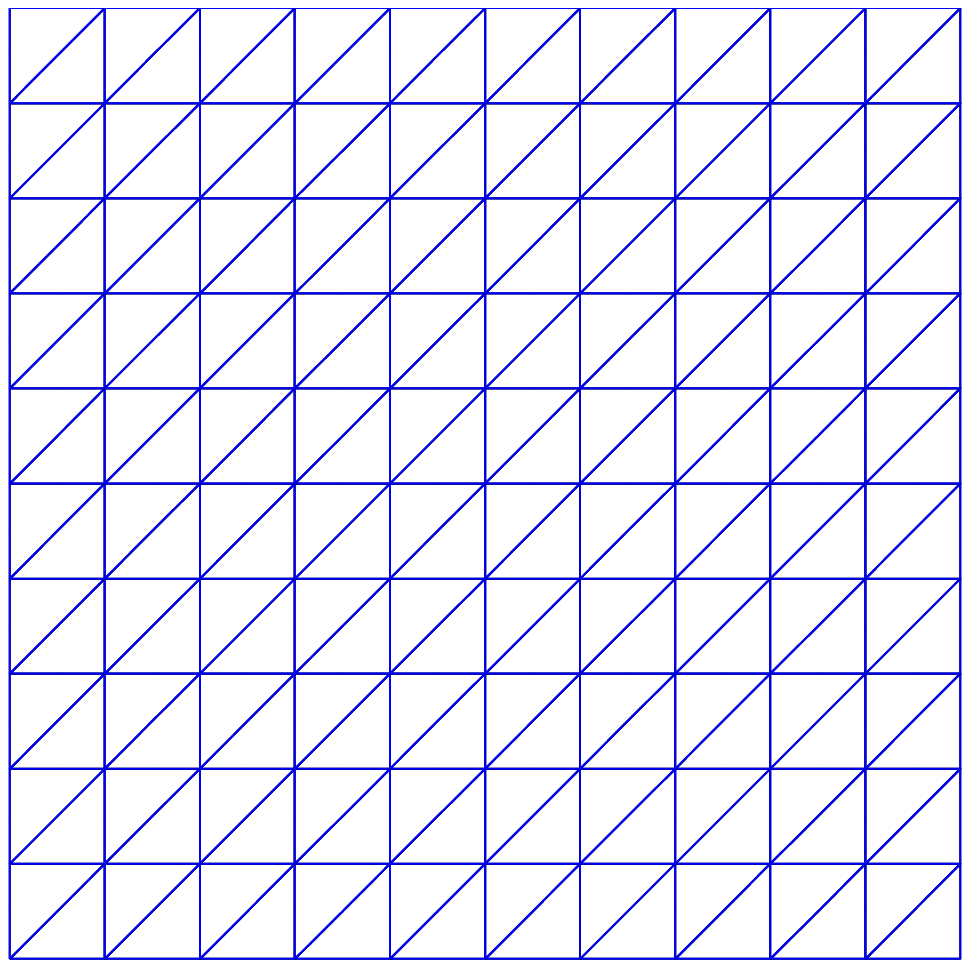}
  \caption{Uniform conforming meshes : mesh 1 (on the left), and mesh 4 (on the right)}
  \label{fig:meshconf}
\end{figure}

\begin{figure}[H]
  \centering
  \hspace{-0.85cm}
   \includegraphics[height=3.85cm]{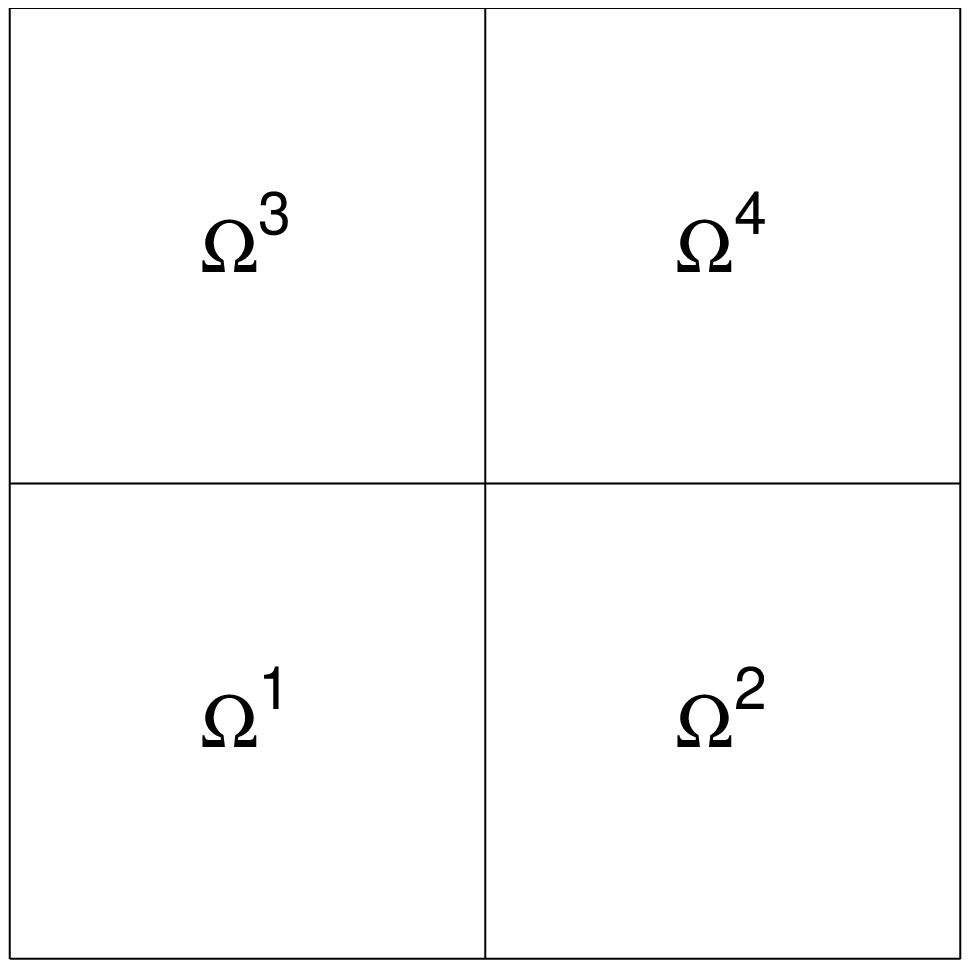}\hspace{-10.mm}
  \includegraphics[height=3.85cm]{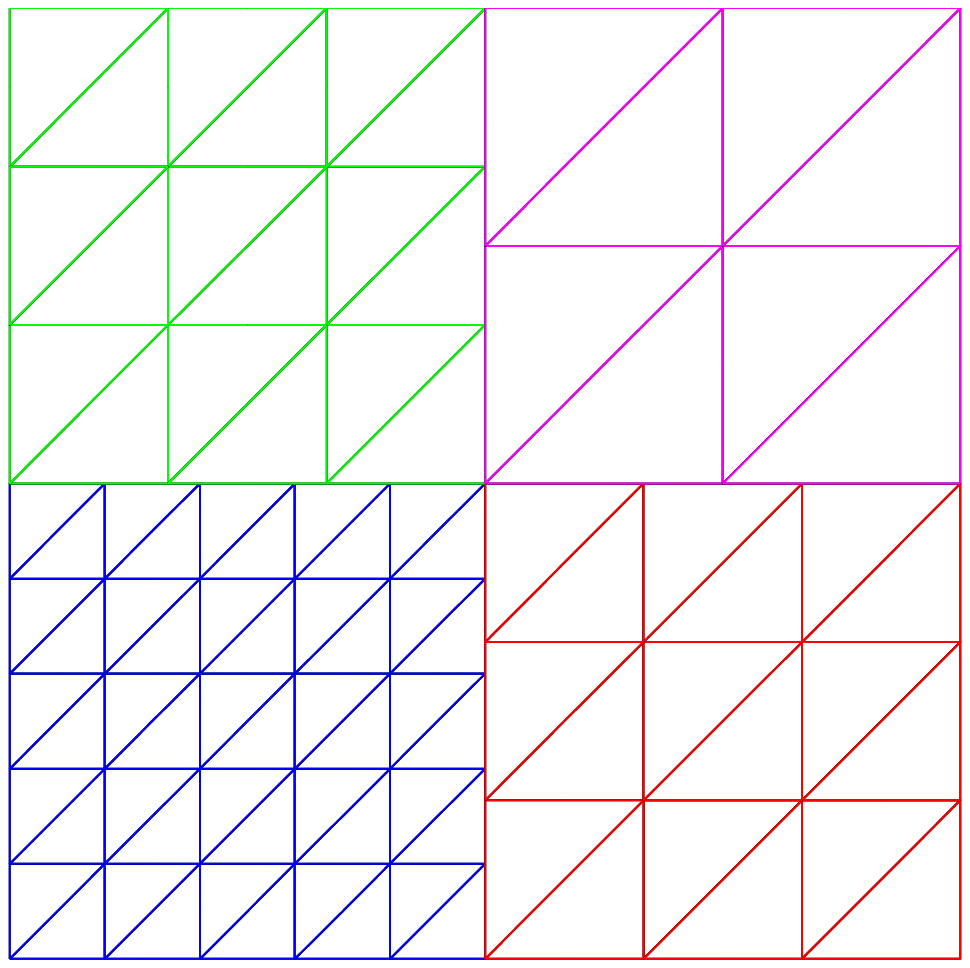}\hspace{-10.0mm}
  \includegraphics[height=3.85cm]{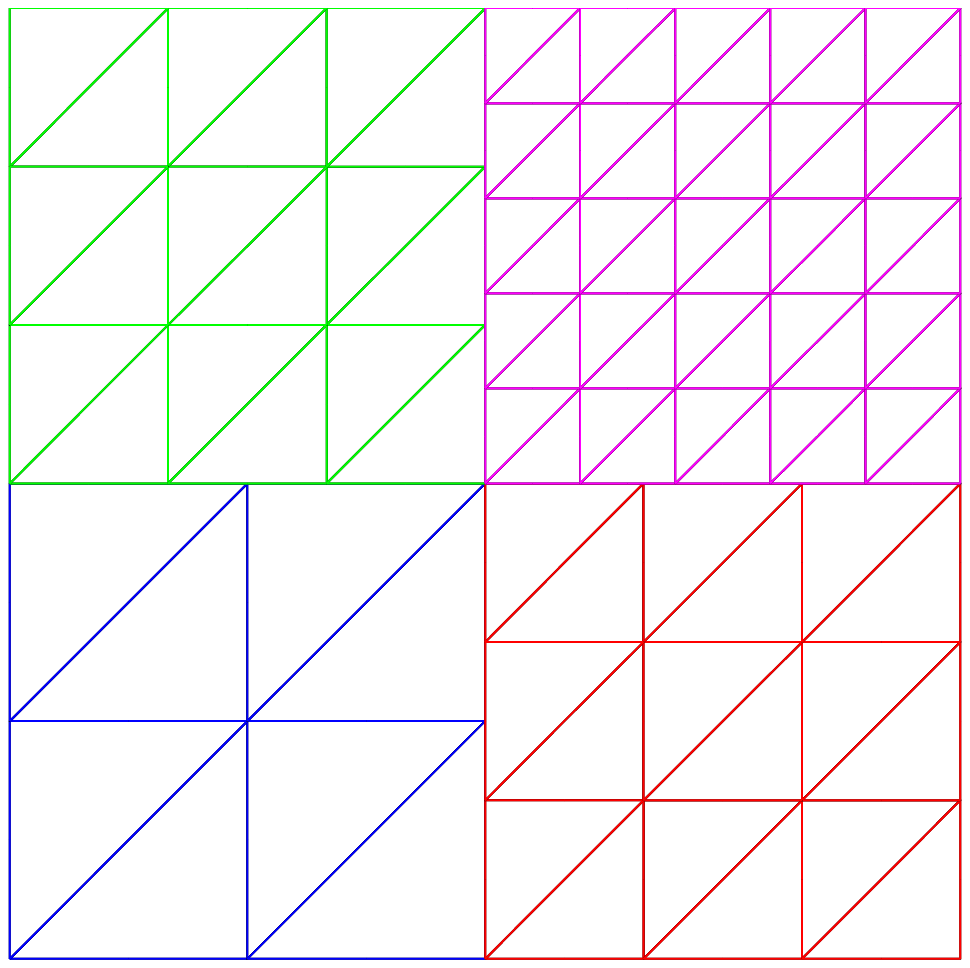}
  \caption{Domain decomposition (on the left), and non-conforming meshes: mesh 2 (on the middle), and mesh 3
    (on the right)}
  \label{fig:meshnc}
\end{figure}

\begin{figure}[H]
  \centering 
  \includegraphics[height=6.5cm]{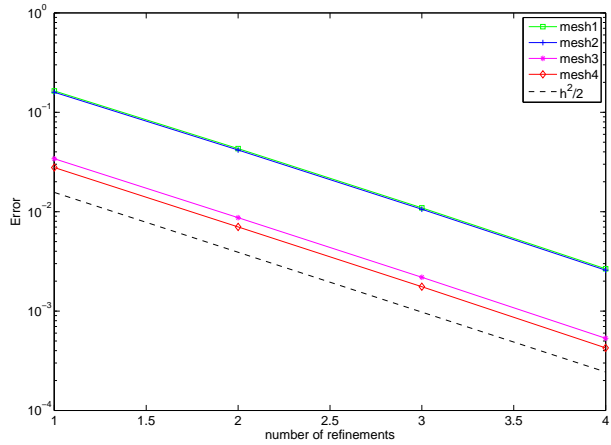}
  \caption{Relative $H^1$ error versus the number of refinements for
    the initial meshes : mesh 1, (square line), mesh 2 (plus line),
    mesh 3 (star line), and mesh 4 (diamond line). The dashed line is
    ${h^2 \over 2}$ (where $h$ is the mesh size) 
    versus the number of refinements, in logarithmic scale}
  \label{fig:errorestim}
\end{figure}

We consider four
initial meshes : the two uniform conforming meshes (mesh 1 and 4) of
Figure~\ref{fig:meshconf}, and the two non-conforming meshes (mesh 2 and 3)
of Figure~\ref{fig:meshnc}.  In the non-conforming case, the unit square is
decomposed into four non-overlapping subdomains numbered as in
Figure~\ref{fig:meshnc} on the left.
Figure~\ref{fig:errorestim} shows the relative $H^1$ error versus the number of
refinement for these four meshes, and ${h^2 \over 2}$ (where $h$ is the mesh size)
versus the number of refinement, in logarithmic scale. At each refinement, the
mesh size is divided by two. The results of Figure~\ref{fig:errorestim} show
that the relative $H^1$ error tends to zero at the same rate as the
mesh size squared~(${h^2}$), and this fits with the theoretical error estimates of
Theorem~\ref{error-estimate2}.

On the other hand, we observe that the two
curves corresponding to the non-conforming meshes (mesh 2 and mesh 3)
are between the curves of the conforming meshes (mesh 1 and mesh
4). We observe that the relative $H^1$ error for mesh 3 is close to the relative $H^1$ error for mesh 4
(i.e. the uniform conforming finer mesh), while the one corresponding to mesh 2 is nearly the same as
the error for mesh 1 (i.e. the uniform conforming coarser mesh), as can be expected,
as mesh 3 is more refined than mesh 2 in subdomain $\Omega^4$, where the solution steeply varies.
%
%
%
\subsection{$H^1$ relative error for different degrees of the finite element approximation}
\label{subsec:errP123}
In this part we study the relative $H^1$ error between the continuous and discrete solutions
versus the mesh size, for $\PP_1$, $\PP_2$ and $\PP_3$ finite elements.
\subsubsection{Decomposition into four subdomains}
\label{subsubsec:err4dom}
\begin{figure}[H]
  \centering
   \hspace{-0.85cm}
   \includegraphics[height=3.85cm]{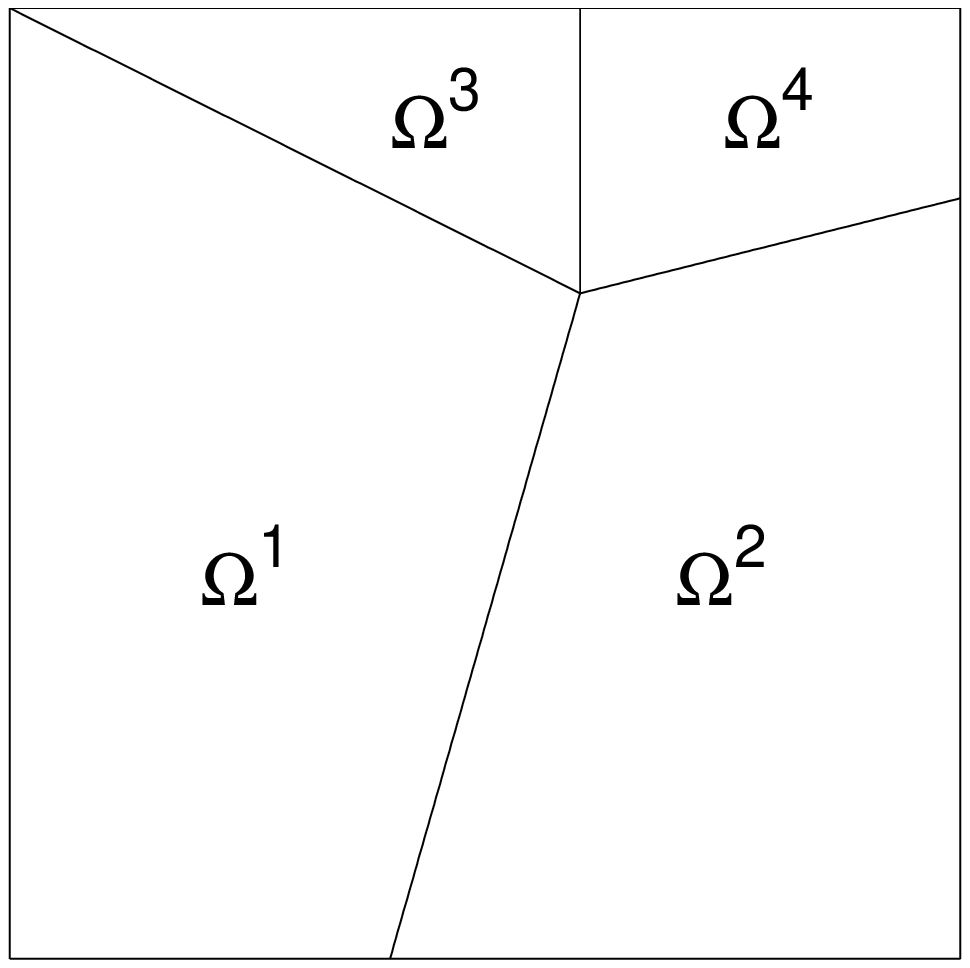}\hspace{-10.mm}
  \includegraphics[height=3.85cm]{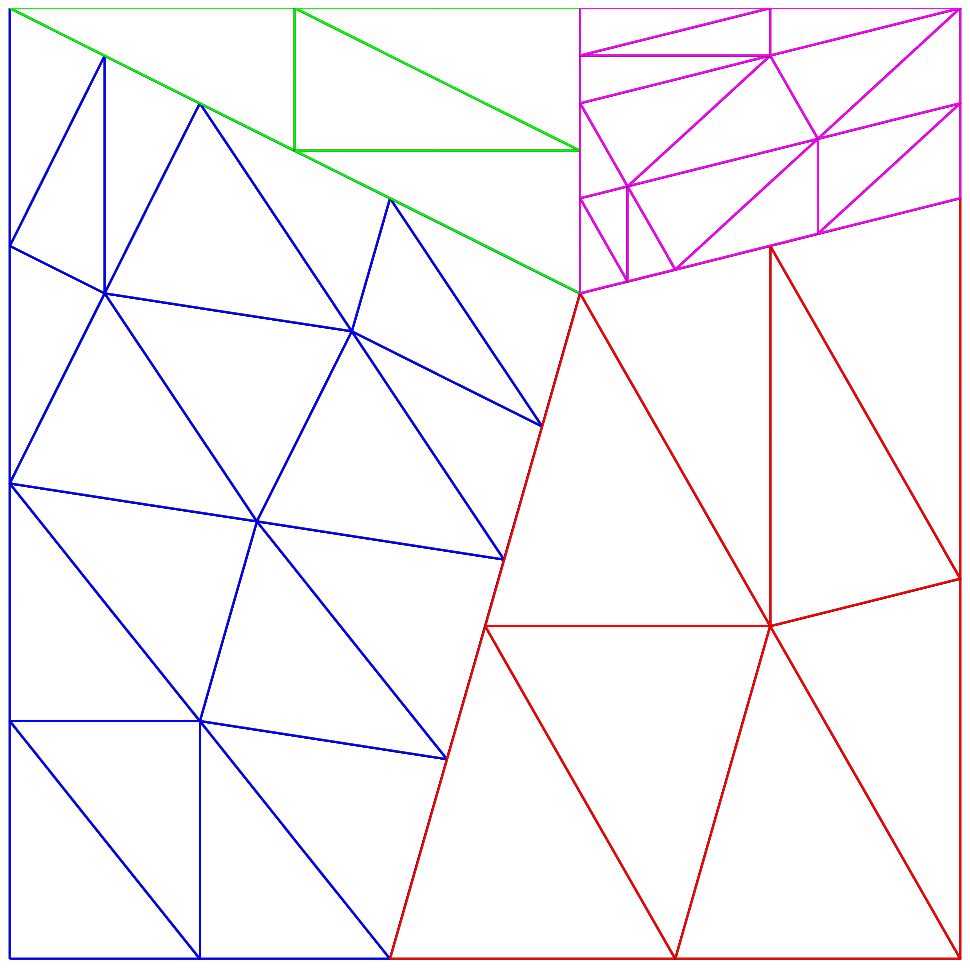}\hspace{-10.mm}
  \includegraphics[height=3.85cm]{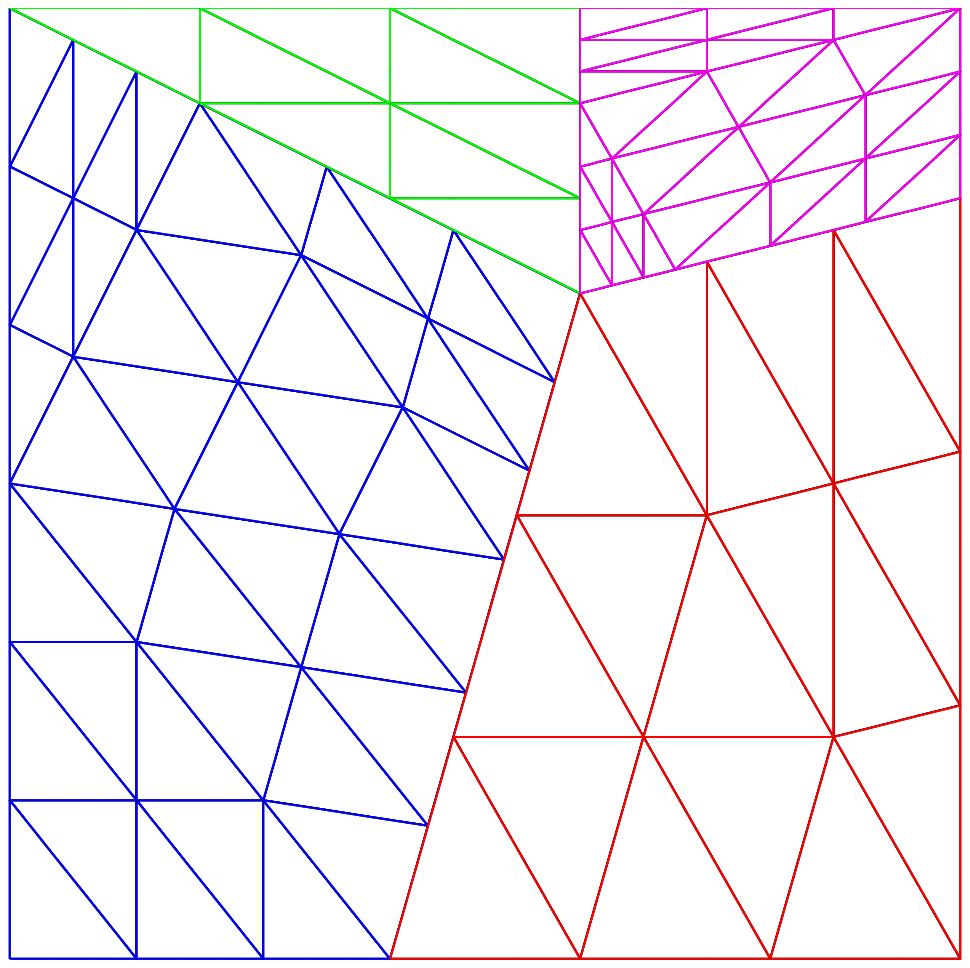}
  \caption{Domain decomposition (left), and non-conforming meshes (middle and right)}
  \label{fig:meshnc4dom}
\end{figure}

\begin{figure}[H]
 \hspace{5mm}
  \includegraphics[height=7.3cm]{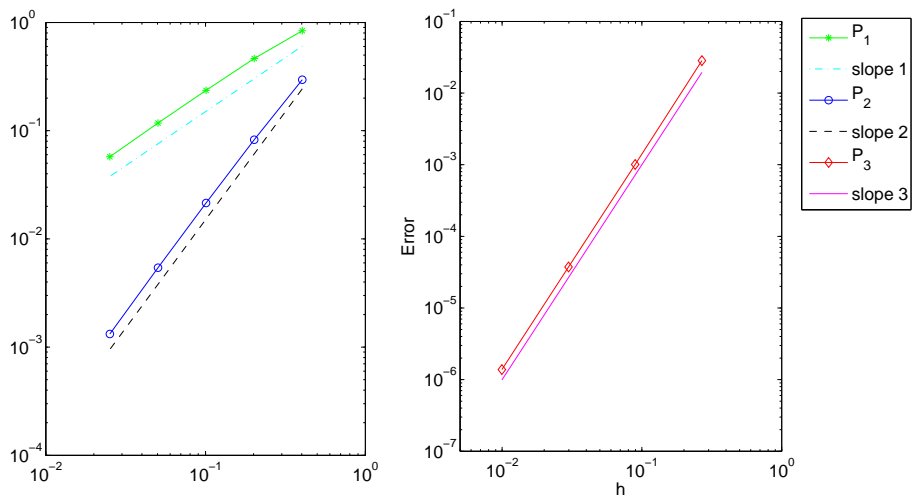}\hspace{-1.5cm}
 \vspace*{-20pt}
  \caption{Relative $H^1$ error versus the mesh size for
     the non-conforming case. Left: for $\PP_1$ and $\PP_2$ discretizations. Right: for $\PP_3$ discretizations}
  \label{fig:errorestim4dom}
\end{figure}

We consider a decomposition of the unit square into four non-overlapping subdomains numbered as in
Figure~\ref{fig:meshnc4dom} on the left.
For $\PP_1$ and $\PP_2$ discretizations, we consider the initial non-conforming meshes represented on
Figure~\ref{fig:meshnc4dom} on the middle, and for a $\PP_3$ discretization, we consider the
initial non-conforming meshes represented on Figure~\ref{fig:meshnc4dom} on the right.
Figure~\ref{fig:errorestim4dom} shows the relative $H^1$ error between the continuous and discrete solutions
versus the mesh size, on the left for $\PP_1$ and $\PP_2$ finite elements, and on the right for
$\PP_3$ finite elements, in logarithmic scales. For $\PP_1$ and $\PP_2$ discretizations,
we start with the meshes on Figure~\ref{fig:meshnc4dom} on the middle
and divide by 2 the mesh size four times. In order to compute the error, the
non-conforming solutions are interpolated on a very fine
grid obtained by refining 5 times the initial mesh. 
For $\PP_3$ discretizations,
we start with the meshes on Figure~\ref{fig:meshnc4dom} on the right
and divide by 3 the mesh size three times. In order to compute the error,
the non-conforming solutions are interpolated on a very fine grid obtained by refining 4 times the initial meshes.

\vspace{-4mm}
The results of Figure~\ref{fig:errorestim4dom} show
that if $p$ is the degree of the approximation, the relative $H^1$ error tends to zero at the same rate as 
$h^p$, for $1\le p \le 3$, and this fits with the theoretical error estimates of
Theorem~\ref{error-estimate2}.

\subsubsection{Decomposition into twelve subdomains}
\label{subsubsec:err12dom}
\begin{figure}[H]
  \centering
  \psfrag{O1}{$ \Omega^1$}
  \psfrag{O2}{$ \Omega^2$}
  \psfrag{O3}{$ \Omega^3$}
  \psfrag{O4}{$ \Omega^4$}
  \psfrag{O5}{$ \Omega^5$}
  \psfrag{O6}{$ \Omega^6$}
  \psfrag{O7}{$ \Omega^7$}
  \psfrag{O8}{$ \Omega^8$}
  \psfrag{O9}{$ \Omega^9$}
  \psfrag{O10}{$ \Omega^{10}$}
  \psfrag{O11}{$ \Omega^{11}$}
  \psfrag{O12}{$ \Omega^{12}$}     
   \includegraphics[height=4.cm]{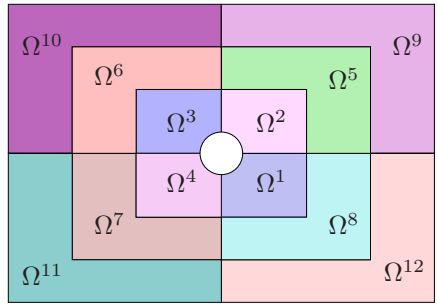}
   \caption{Domain decomposition into twelve non-overlapping subdomains}
  \label{fig:dec12dom}
\end{figure}

\begin{figure}[H]
  \centering
   \includegraphics[height=8cm]{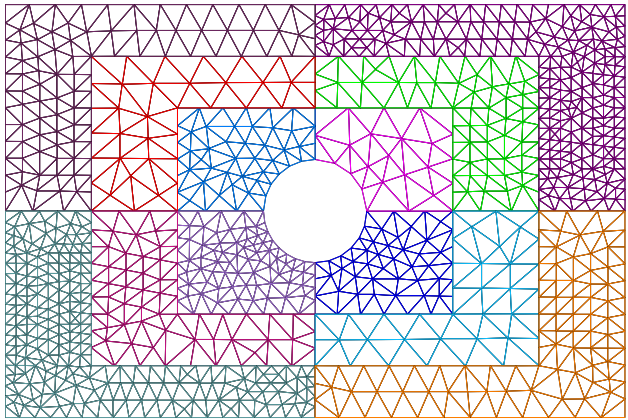}
   \vspace*{-30pt}
  \caption{Non-conforming meshes}
  \label{fig:meshnc12dom}
\end{figure}

\begin{figure}[H]
 \hspace{-25mm}
  \includegraphics[height=8.5cm]{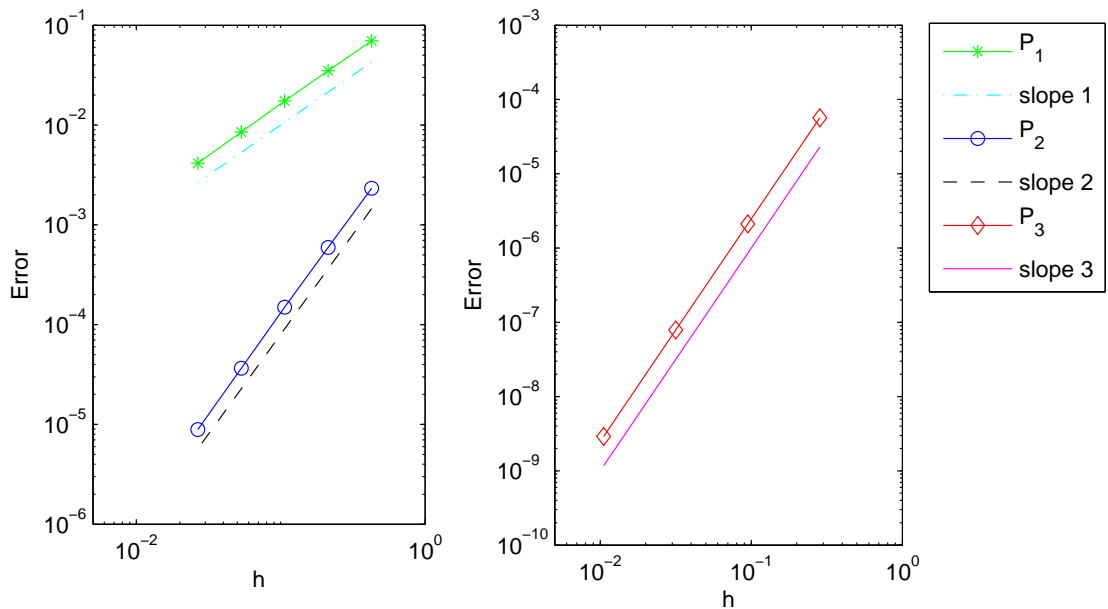}\hspace{-1.5cm}
 \vspace*{-20pt}
  \caption{Relative $H^1$ error versus the mesh size for
     the non-conforming case. Left: for $\PP_1$ and $\PP_2$ discretizations. Right: for $\PP_3$ discretizations}
  \label{fig:errorestim12dom}
\end{figure}

We consider the initial problem with  exact solution $u(x,y)=x^3y^2+\sin(xy)$.
The domain is $\Omega=(-3,3) \times (-2,2)$ and is decomposed into
twelve irregularly shaped subdomains as in Figure~\ref{fig:dec12dom}.
The subdomain meshes are generated in an
independent manner as in Figure~\ref{fig:meshnc12dom}. The finite element assemblies are done as in \cite{Cuvelier}.
\vspace*{-5mm}
Figure~\ref{fig:errorestim12dom} shows the relative $H^1$ error between the continuous and discrete solutions
versus the mesh size, on the left for $\PP_1$ and $\PP_2$ finite elements, and on the right for
$\PP_3$ finite elements, in logarithmic scales.

For $\PP_1$ and $\PP_2$ discretizations,
we start with the mesh on Figure~\ref{fig:meshnc12dom} and divide by 2 the mesh size four times.
In order to compute the error, the non-conforming solutions are interpolated on a very fine
grid obtained by refining 5 times the initial mesh.
For $\PP_3$ discretizations, we start with the mesh on Figure~\ref{fig:meshnc12dom} 
and divide by 3 the mesh size three times.
In order to compute the error, the non-conforming solutions are interpolated on a very fine grid obtained
by refining 4 times the initial meshes.

The results of Figure~\ref{fig:errorestim12dom} show
that the relative $H^1$ error tends to zero at the same rate as 
$h^p$, for $1\le p \le 3$ where $p$ is the degree of the approximation. This corresponds to the theoretical
error estimates of Theorem~\ref{error-estimate2}.
\vspace*{-5mm}
%
%
\subsection{Convergence : Choice of the Robin parameter}\label{sec.Robinparam}
\begin{figure}[H]
  \centering
  \includegraphics[height=4.cm]{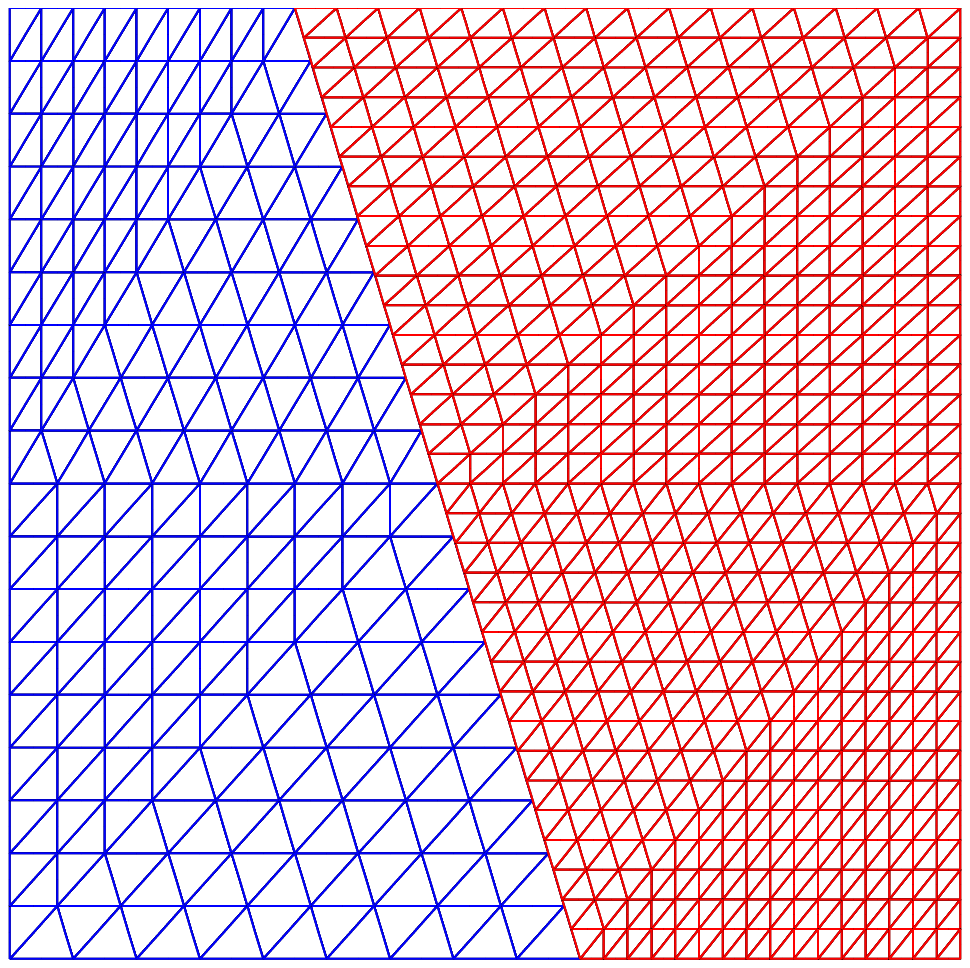}
  \includegraphics[height=4.cm]{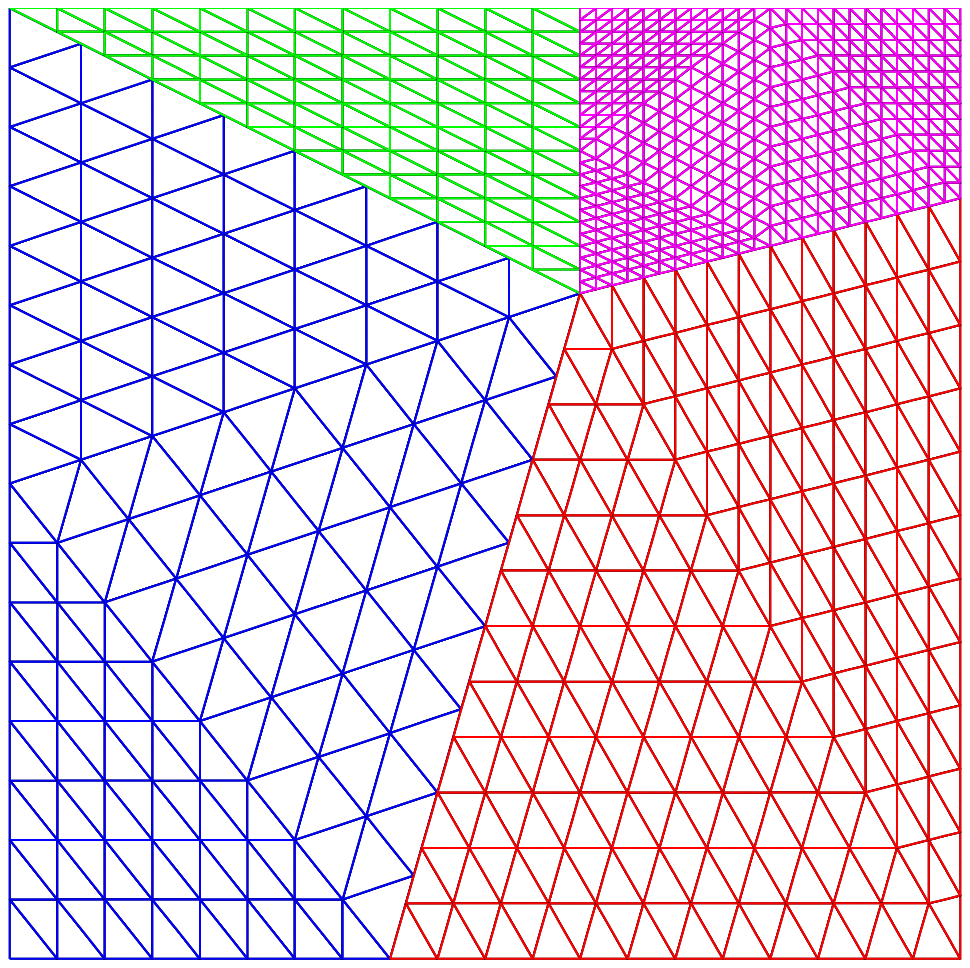}
  \caption{Domain decomposition in two subdomains (left) and in four subdomains (right), with non-conforming meshes}
  \label{fig:meshnc2-4dom}
\end{figure}

%
Let us now study the convergence speed to reach the discrete solution, for
different values of the Robin parameter $\alpha$, in the case of $\PP_2$ finite elements.
We first consider a domain decomposition in two subdomains, and then
in four subdomains as shown in Figure~\ref{fig:meshnc2-4dom}. We simulate the error equations (i.e. $f = 0$),
and use a random initial guess so that all the frequency components are present.
\subsubsection{2 subdomain case}\label{sec.2sub}
In this part, the unit square is decomposed into two subdomains
with non-conforming meshes (with $703$ and $2145$ nodes respectively)
 as shown in Figure~\ref{fig:meshnc2-4dom} (on the left). 

\begin{figure}[H]
  \hspace*{-1.8cm}
  \includegraphics[height=6.4cm]{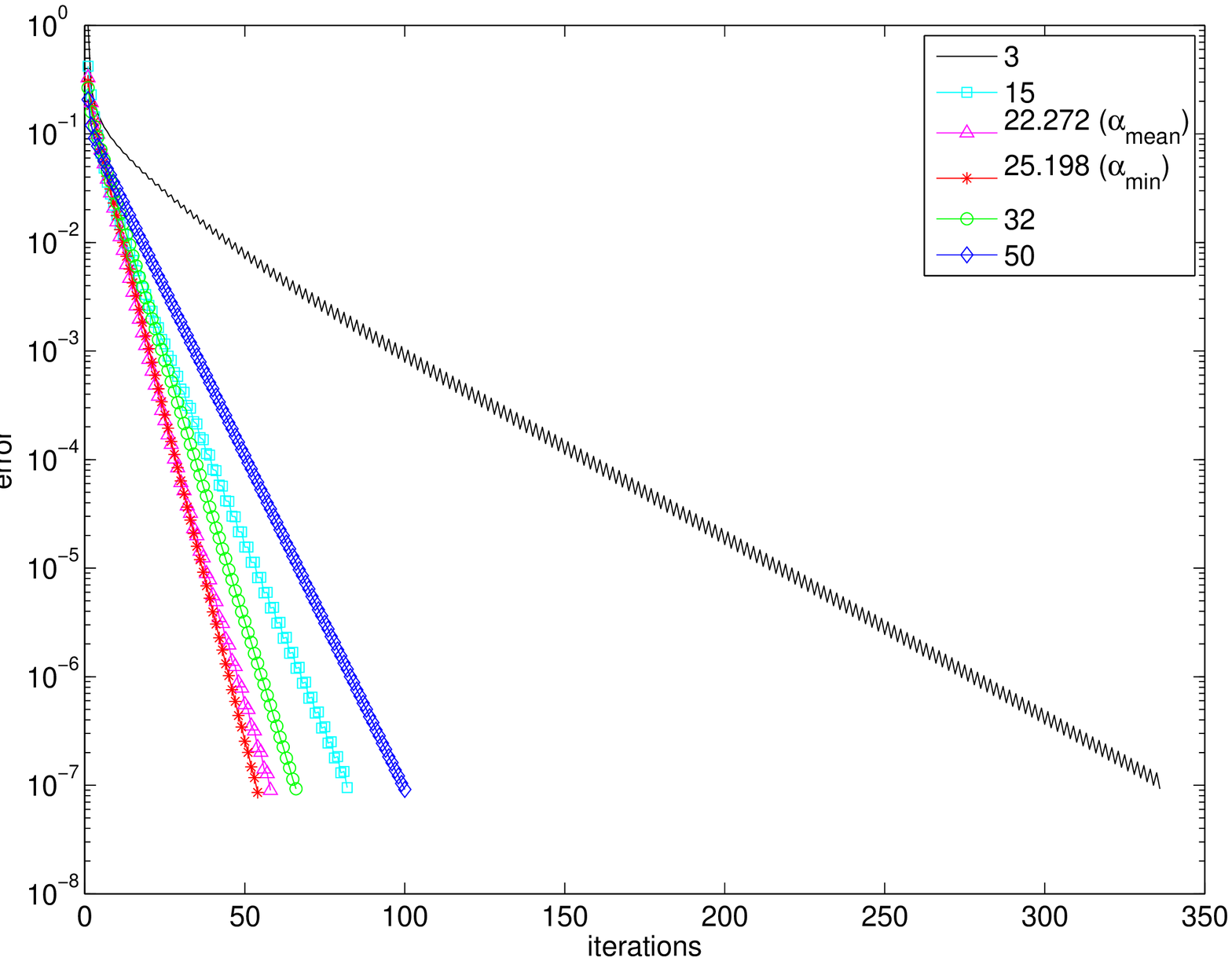}\hspace{-0.7cm}
    \includegraphics[height=6.4cm]{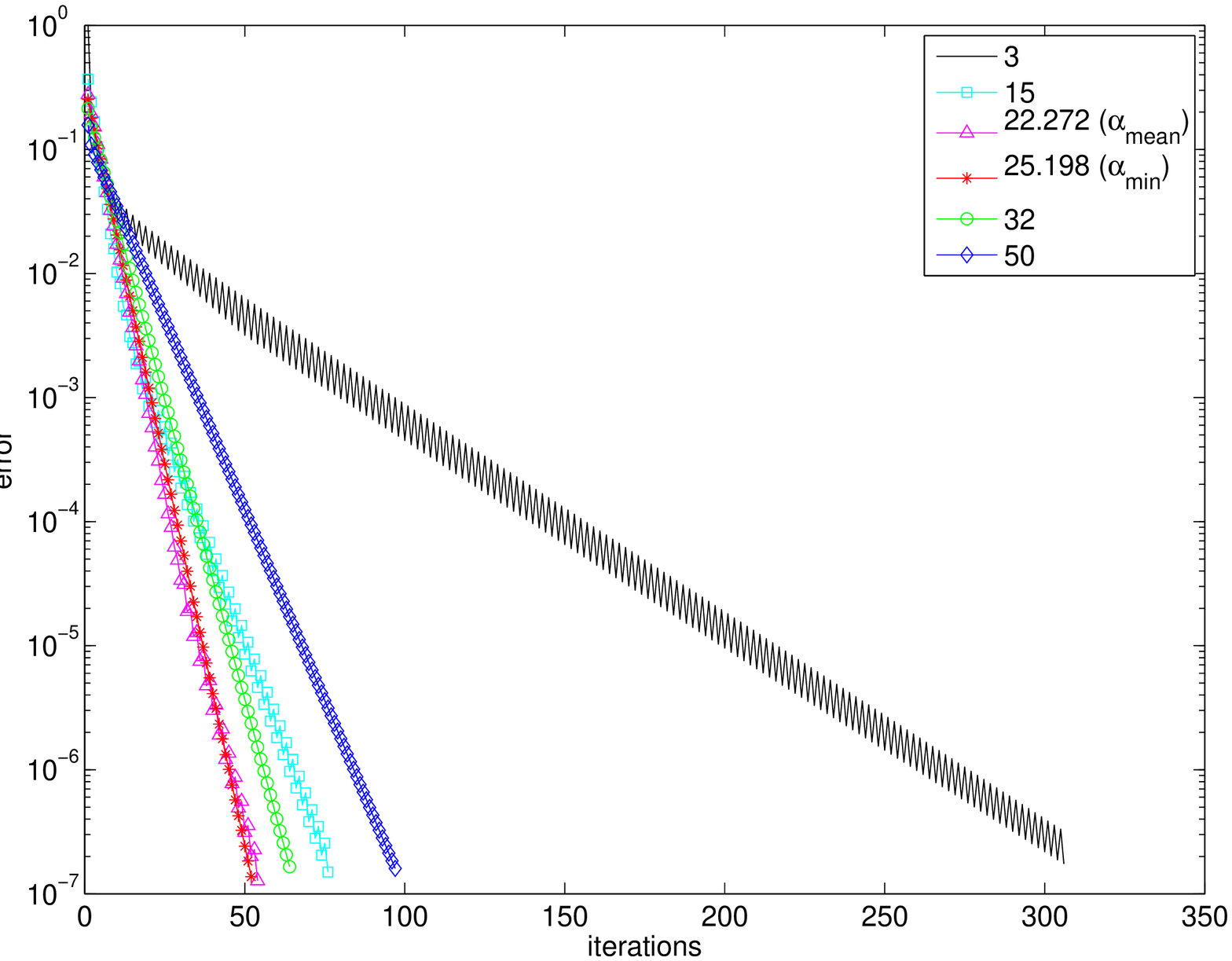}
    \hspace*{-1.8cm}
 \includegraphics[height=6.4cm]{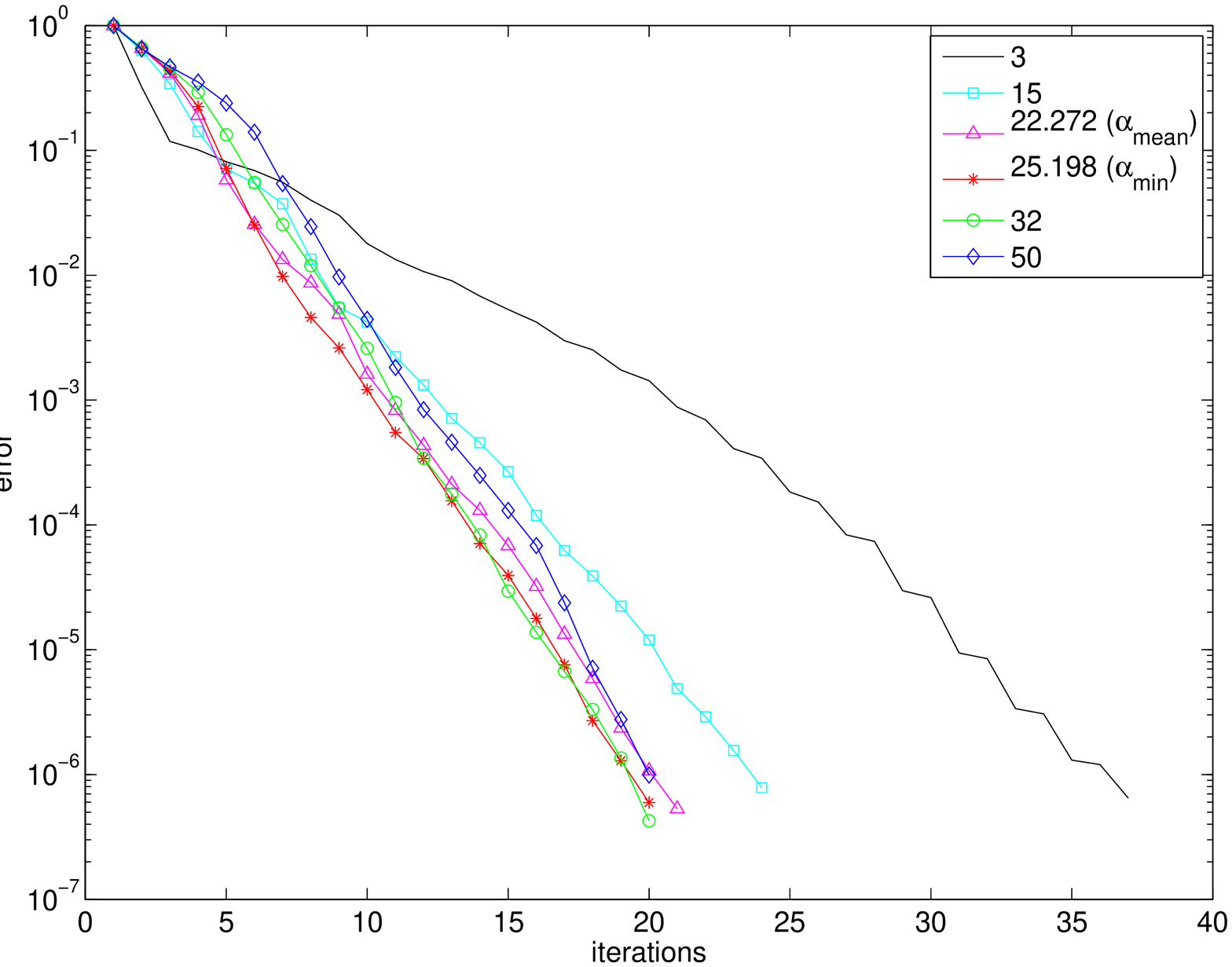}\hspace{-0.7cm}
  \includegraphics[height=6.4cm]{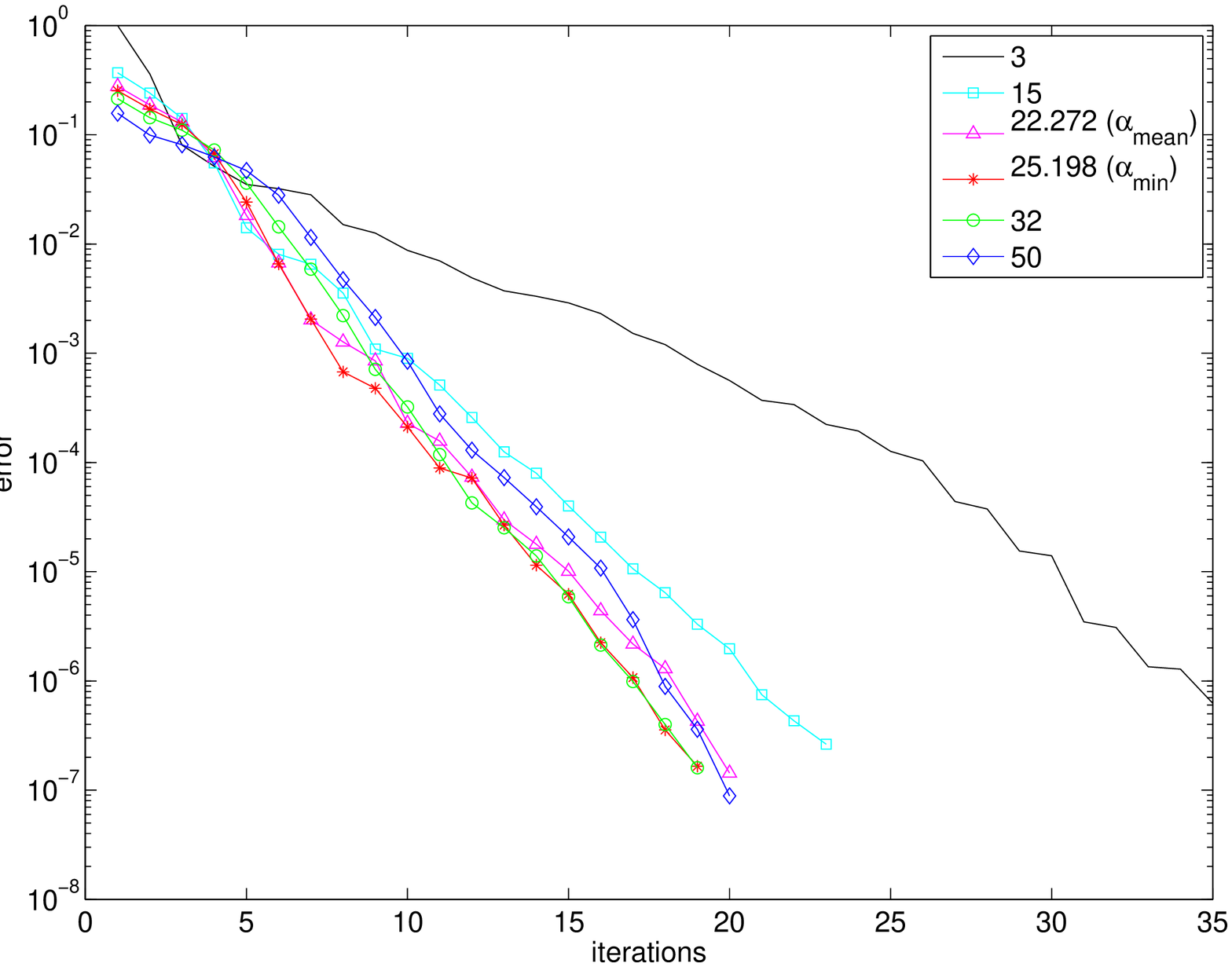} 
 \caption{Error versus Schwarz (top) or Gmres (bottom) iterations for different values of the Robin 
parameter $\alpha$. Left: the $H^1$ error, right: the $L^\infty$ error.}
\label{fig:errconv2dom}  
\end{figure}

On Figure~\ref{fig:errconv2dom} (top left) we represent
 the $H^1$ norm of the iterate error, for different values of the Robin 
parameter $\alpha$. We observe that
the optimal numerical value of the Robin parameter is close to $\alpha_{\text{min}}$.
As the relative $H^1$ error didn't show where the error is highest,
we also represented on Figure~\ref{fig:errconv2dom} (top right)
the $L^{\infty}$ norm of the iterate error, for different values of the Robin parameter $\alpha$.
We obtain similar results as for the relative $H^1$ error.

The Schwarz algorithm can be interpreted as a Jacobi algorithm applied
to an interface problem (see~\cite{Nataf.4}). In order to accelerate the 
convergence, we can replace the Jacobi algorithm by a GMRES~(\cite{Saad}) algorithm.
Figure~\ref{fig:errconv2dom} show respectively the
$H^1$ norm (on the bottom left) and the $L^{\infty}$ norm (on the bottom right) of the GMRES iterate error,
for different values of the Robin parameter $\alpha$.
\vspace{-3mm}
For $\alpha=\alpha_{\text{min}}$, the convergence is 
accelerated by a factor 2 for GMRES, compared to the Schwarz algorithm.
Moreover, the gap between the error values for different $\alpha$  
is decreasing when using the GMRES algorithm, compared to the Schwarz method.
Thus, the GMRES algorithm is less sensitive to the choice of the Robin
parameter. The sensitivity of the performance of the Krylov solver to the
optimized value of the parameter is thus not so critical but it is
real and especially visible for ranges of accuracy used for most
practical applications (relative errors of size $10^{-2}$ or
$10^{-3}$). Adding that this effect generally increases with the number
of subdomains and the refinement of the mesh \cite{Gander06} together with the complexity of the equations,
we advise when possible to look for the optimized value.
Moreover, this conclusion on the interest of GMRES compared to Schwarz is established for stationnary problems but is
not yet verified for time dependent problems with a Schwarz waveform relaxation algorithm,
as illustrated for example in~\cite{Hoang}.

In Table~\ref{table:errconv} we show the number of iterations $N$ to reduce the $H^1$ error by a factor
  $10^6$ versus the Robin parameter $\alpha$, for different degrees $p$ of the approximation. We observe that
  $\alpha_{\min}$ is very close to the optimal numerical value, for all $p=1,2,3$.
%
\begin{table}[H]
\begin{center}
\begin{tabular}{|l|l||l|l||l|l|}
\hline
$\alpha$ & N ($p=1)$ & $\alpha$ & N ($p=2$)  & $\alpha$ & N ($p=3$)\\ \hline \hline
10 & 57 & 17 & 63 & 23 & 88 \\ \hline 
15 & 40 & 22 & 51 & 28 & 74\\ \hline
$17.818 \,(\alpha_{\min})$ & 36 & $25.198 \,(\alpha_{\min})$ & 49 & $30.861 \,(\alpha_{\min})$ & 68 \\ \hline
20 & 39 & 27 & 52 & 33 & 66\\ \hline
25 & 40 & 32 & 61 & 35& 68\\ \hline
30 & 58 & 37 & 71 & 40& 77\\ \hline
\end{tabular}
\end{center}
\caption{Number of iterations to reduce the $H^1$ error by a factor
  $10^6$ versus $\alpha$, for different degrees $p$}
\label{table:errconv}
\end{table}

\subsubsection{4 subdomain case}
In this part, the unit square is decomposed into four subdomains
with non-conforming meshes as shown in Figure~\ref{fig:meshnc2-4dom} (on the right).

\begin{figure}[H]
  \hspace{-1.5cm}
  \includegraphics[height=6.cm]{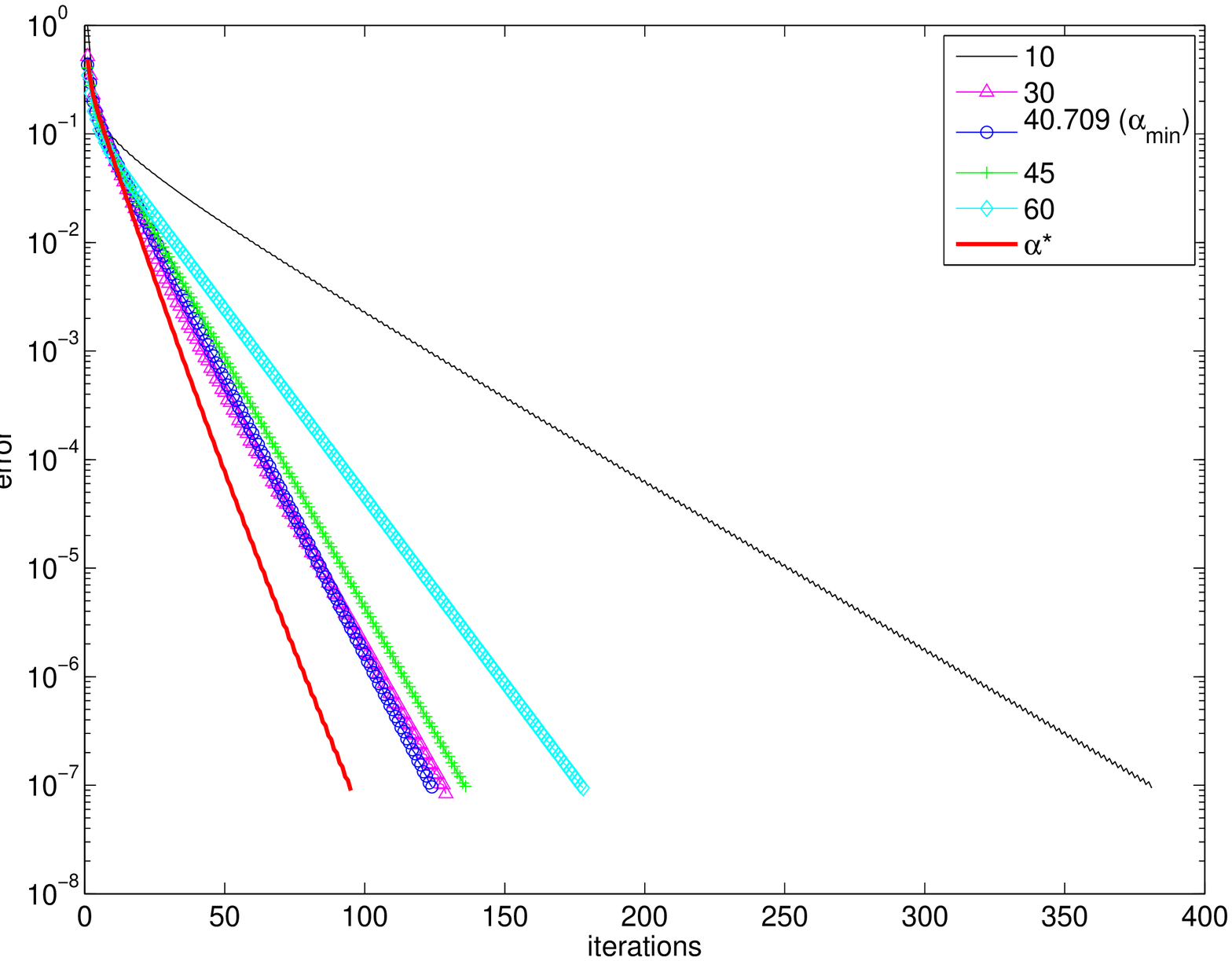}\hspace{-0.7cm}
    \includegraphics[height=6.cm]{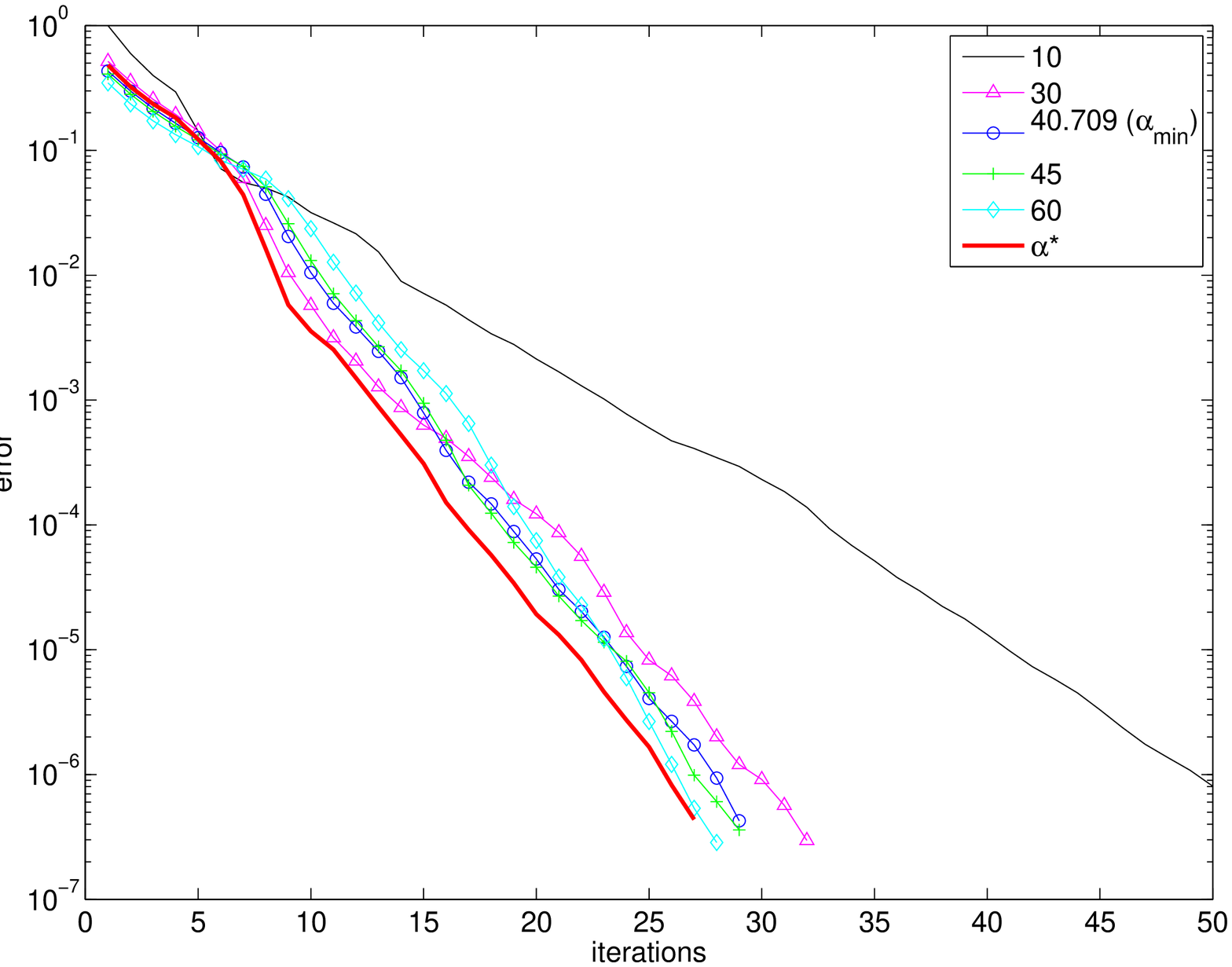}
    \caption{Error versus iterations, for different values of the Robin 
parameter $\alpha$. Left: the Schwarz algorithm, right: the GMRES algorithm.}
\label{fig:err4d}       
\end{figure}

From the results of Section~\ref{sec.2sub}, we will consider for the optimized parameter $\alpha$ the values
given by the smallest mesh size on the interface.
As we have four interfaces, using formula \eqref{eq.alphaopt}
with $h=h_{\text{min}}^{k,\ell}$, $ 1\le k,\ell \le 4$, with $\Gamma^{k,\ell}$
not empty, we obtain four values given by:
$\alpha_{\text{min}}^{1,2}=18.479,\ \alpha_{\text{min}}^{1,3}=19.250, \ \alpha_{\text{min}}^{2,4}=34.725,
\ \alpha_{\text{min}}^{3,4}=40.709$.
We define $\alpha^*$ the parameter with these four values over the interfaces (i.e. $\alpha^*$ is constant
over each interface, with different constants from one interface to another).
We consider also a constant optimized value $\alpha_{\text{min}}$ over the four interfaces obtained 
by taking $h=\min(h_{\text{min}}^{1,2},h_{\text{min}}^{1,3},h_{\text{min}}^{2,4},h_{\text{min}}^{3,4})$ in
formula \eqref{eq.alphaopt}. We obtain $\alpha_{\text{min}}=\alpha_{\text{min}}^{3,4}=40.709$.
 On Figure~\ref{fig:err4d} we represent
the $H^1$ norm of the error, for $\alpha^*$, $\alpha_{\text{min}}$ and for different constant values of the Robin 
parameter $\alpha$, with the Schwarz method on the left, and the GMRES algorithm on the right.

We observe that the optimal numerical value of the Robin parameter is close to $\alpha_{\text{min}}$.
We also observe in that case that taking the different optimized values over the interfaces (i.e. $\alpha^*$)
doesn't improve substantially the convergence speed compared to taking the same value $\alpha_{\text{min}}$ over
the interfaces.
\section{Conclusions}
%
We have analyzed the convergence of the iterative algorithm
for $\PP_p$ finite elements, with $p\ge 1$ in 2D and $p=1$ in 3D,
for the NICEM method. It relies on Schwarz type algorithms with Robin
interface conditions on non-conforming grids. We have extended the error estimates in 2D for piecewise polynomials
of higher order. Numerical results show that the method preserves the order of the finite elements
for $\PP_p$ discretizations, with $p=1,2$ or $3$.

\appendix
\section{Proof of Lemma~\ref{lem:Delta} in the case $\pmb{p\ge 3}$}\label{appendix:A}
From Section~\ref{sec:peq2}, it remains 
  to prove Lemma~\ref{lem:Delta} in the general case $p\ge 3$.
Let us introduce the vector space $Q^p=\{
\eta\in\PP_p([-1,1])
\hbox{ s.t. }\
\eta(-1)=0\}$.  The function $\Delta(\eta)$ is quadratic so that it
suffices to study the extrema of $\Delta(\eta)/\|\eta\|^2_{L^2(]-1,1[)}$
over $Q^p$ or equivalently to prove that the associated symmetric 
quadratic form is
negative, i.e. its eigenvalues are negative. They correspond to the Lagrange
multiplier solutions
$\mu_1$ of the following min-max problem
\begin{equation}\label{eq:minMaxDiscr}
\min_{\eta\in Q^p} \max_{\mu_1\in\R} {\cal L}_e(\eta,\mu_1),
\end{equation}
where
\[
{\cal L}_e(\eta,\mu_1) := \Delta(\eta) - \mu_1 (\|\eta\|_{L^2(]-1,1[)}^2-1).
\]
We have to prove that $\mu_1<0$. We have 
\[\hspace{-0.3cm}\begin{array}{rcl}
0&=&\ds<\frac{\partial{\cal L}_e}{\partial \eta},\delta \eta >\\
\ds\hphantom{ddd}&=&
2(2\eta(1)-3\eta_p)(2\delta\eta(1)-3\delta\eta_p)+p^2(-8<\eta,\delta\eta>
+18\frac{\eta_p\delta\eta_p}{2p+1})
-2\mu_1 <\eta,\delta\eta>
\end{array}\]
where $<\, ,\,>$ denotes the $L^2$ scalar product on $L^2(]-1,1[)$ and
$\delta\eta\in Q^p$.
\\ Let us consider the vector space $(1-x^2) \PP_{p-3}\subset Q^p$. Any
  function $\gamma$ in $(1-x^2) \PP_{p-3}$ satisfies $\gamma(-1)=\gamma(1)=0$
and $\gamma_p=0$.
The optimality relation w.r.t. to $(1-x^2) \PP_{p-3}$ gives
\[
(-8p^2-2\mu_1)<\eta,\delta\eta>=0,\ \ \forall
\delta\eta\in (1-x^2) \PP_{p-3}.
\]
We have either $\mu_1=-4p^2<0$ or  $\eta$  solution to (\ref{eq:minMaxDiscr})
belongs to the space $\{(1-x^2) \PP_{p-3}\}^\bot
\cap \PP_p$. The first case corresponds to a negative value for 
$\mu_1$ which is in
agreement with the lemma to be proved. Let us study the latter case. 
We shall make
use of (see~\cite{AS})
\begin{lem}\label{lem:LagrangePrime} The family of Legendre polynomials satisfies
\begin{eqnarray*}
  \begin{array}{r}
\ds\int_{-1}^1 L'_m\,L'_{m'}\,(1-x^2) \, dx = 0, \ m\neq m',\\
\ds\int_{-1}^1 {L'_m}^2= m(m+1),\\
\ds\int_{-1}^1 L'_m\,L'_{m+1}= 0,\\
\ds\int_{-1}^1 L'_{m-1}\,L'_{m+1}= m(m-1),\\
\ds L'_m(-1) = (-1)^{m+1} \frac{m(m+1)}{2}.
\end{array}
\end{eqnarray*}
For any $p$, $p\ge 3$,
\[
\{(1-x^2) \PP_{p-3}\}^\bot \cap \PP_p = \mathrm{ Span}\{L_{p},\,L'_p,\,L'_{p-1}\}.
\]
\end{lem}
{\bf Proof.} We only need to prove the last equality, that results from the above indeed it can be checked easily that
\[
\{(1-x^2) \PP_{p-3}\}^\bot \cap \PP_p = \hbox{ Span}\{L'_{p+1},\,L'_p,\,L'_{p-1}\}.
\]
Moreover, we have
\[
L'_{p+1}(x)=(2p+1)L_p(x)+L'_{p-1}(x)
\]
and thus Lemma~\ref{lem:LagrangePrime}.\\
Therefore, there exists $\lambda_1,\lambda_2,\lambda_3 \in\R$ s.t.
$\eta=\lambda_1 L_{p} + \lambda_2 L'_{p} + \lambda_3 L'_{p-1}$.
Since $\eta$ is defined up to a constant, we only
have to consider the two cases
$\lambda_1=1$ or $\lambda_1=0$.\\[1.5mm]
{\bf Case 1:} $\lambda_1=1$\\
 From $\eta(-1)=0$, we get
\[
1-\lambda_2 \frac{p(p+1)}{2} +\lambda_3 \frac{p(p-1)}{2} =0,
\]
so that
\[
\lambda_2 = \frac{2}{p(p+1)} +\lambda_3  \frac{p-1}{p+1},
\]
\begin{eqnarray*}
\Delta(\eta)=\displaystyle - 4 \frac{(p-1)\,p^{2}\,(p^{2} +  1)}{p + 1}\lambda_3^{2}  -  
\frac {(24\,p^{4} - 20\,p^{3} - 8\,p^{2} + 4\,p)}{(p + 1)\,(2\,p + 1)}\lambda_3 
-\frac {29\,p^{2} + 13\,p - 1 - p^{3}}{(p + 1)\,(2\,p + 1)}.
\end{eqnarray*}
Since $p$ is supposed larger than~1, the
leading coefficient of $\Delta(\eta)$ is negative. If the discriminant of
$\Delta(\eta)$ is negative, the polynomial is negative for any
$\lambda_3$. This discriminant has the value
\[
16\,{\displaystyle \frac {(p^{2} - 13\,p - 8)\,(p - 1)\,p^{3}}{2
\,p + 1}}
\]
and is negative for $2\le p \le 13$.\\[1.5mm]
{\bf Case 2:} $\lambda_1=0$\\
 From $\eta(-1)=0$, we get
\[
-\lambda_2 \frac{p(p+1)}{2} +\lambda_3 \frac{p(p-1)}{2} =0,
\]
so that
\[
\lambda_2 = \lambda_3  \frac{p-1}{p+1}.
\]
Since $\eta$ is an eigenvalue, it is not zero and the above relation
shows that we can take $\lambda_3=1$. Then, we have
$\lambda_2 = \frac{p-1}{p+1}$ so that
$$
\Delta(\eta)=\frac{-4(p-1)p^2(p^2+1)}{(p+1)} < 0,
$$
which ends the proof of Lemma~\ref{lem:Delta}. $\qquad \Box$

\section{Proof of the estimate \pmb{\eqref{eq:logh}}}\label{appendix:B}
For $0\le m < p-1, \ \beta(m)=0$ and the estimate (\ref{eq:logh}) is standard.
For $m=p-1$, the proof is the same as for Lemma~5 in~\cite{JMN10}: let ${\bar p}_{k \ell h}$ be the unique element
of $\tilde W_h^{k,\ell}$ defined as follows :
\begin{itemize}
\item
$({\bar p}_{k \ell h})_{|[x_1^{\ell,k},x_{N-1}^{\ell,k}]}$ coincide
with the interpolate of degree $p$ of $p_{k,\ell}$.
\item
$({\bar p}_{k \ell h})_{|[x_0^{\ell,k},x_{1}^{\ell,k}]}$ and 
$({\bar p}_{k \ell h})_{|[x_{N-1}^{\ell,k},x_{N}^{\ell,k}]}$ coincide
with the interpolate of degree $p-1$ of $p_{k,\ell}$.
\end{itemize}
Then, using Deny-Lions theorem we have
\begin{align*}
\| p^1_{k \ell h} - p_{k,\ell} \|^2_{L^2((\Gamma^{k,\ell})}
\le& \| {\bar p}_{k \ell h} - p_{k,\ell} \|^2_{L^2(]x_0^{\ell,k},x_1^{\ell,k}[)} \\
&+ch^{1+2p}\|p_{k,\ell}\|^2_{H^{{1\over 2}+p}([x_1^{\ell,k},x_{N-1}^{\ell,k}])}
+\| {\bar p}_{k \ell h} - p_{k,\ell} \|^2_{L^2(]x_{N-1}^{\ell,k},x_N^{\ell,k}[)}.
\nonumber
\end{align*}
In order to analyze the two extreme contributions, we use Deny-Lions theorem 
\bea
\| {\bar p}_{k \ell h} - p_{k,\ell} \|^2_{L^2(]x_0^{\ell,k},x_1^{\ell,k}[)}
\le ch^{1+2p-{2\over q}}\|\frac{d^p p_{k,\ell}}{dx^p}\|^2_{L^{q}(]x_0^{\ell,k},x_1^{\ell,k}[)},
\nonumber
\eea

and taking  $q=-log(h)$, we finish the proof as for Lemma~5 in~\cite{JMN10}.
$\qquad \Box$

\subsection*{Acknowledgment}
The authors would like to thank Fran\c{c}ois Cuvelier for his help in the implementation of the test case of
  Section~\ref{subsubsec:err12dom}, especially
for his development of a FreeFem++ code that generates automatically the meshes with different refinement levels,
that we used for our numerical results.

\end{document}